\documentclass[12pt]{amsart}
\usepackage{latexsym, graphics, amscd, amssymb, amsfonts,
amsmath}
\def\fuchs{\wedge^{\infty}_{}V_{\infty}}
\def\ni{\noindent}

\def\v{\varphi}

\newcommand{\lam}{\lambda}
\newcommand{\Lam}{\Lambda}

\newcommand{\PP}{\mathbb{P}}

\newcommand{\cb}{\mathcal{B}}
\newcommand{\cg}{\mathcal{G}}

\newcommand{\cp}{\mathcal{P}}

\newcommand{\CC}{\mathbb{C}}

\newcommand{\ZZ}{{\mathbb Z}}
\newcommand{\NN}{{\mathbb N}}

\newcommand{\II}{\mathcal{I}}

\newcommand{\FF}{\mathcal{F}}
\newcommand{\HH}{\mathcal{H}}

\newcommand{\NNN}{\mathcal{N}}

\newcommand{\X}{\mathop{\rm X}\nolimits}
\newcommand{\SLhat}{\mathop{\widehat{\rm SL}}\nolimits}

\newcommand{\Grhat}{\mathop{\widehat{\rm Gr}}\nolimits}
\newcommand{\Gr}{\mathop{\rm Gr}\nolimits}
\newcommand{\What}{\mathop{\widetilde{W}}\nolimits}

\newcommand{\ord}{\mathop{\rm ord}\nolimits}

\newcommand{\vdim}{\mathop{\rm vdim}\nolimits}
\newcommand{\vcard}{\mathop{\rm vcard}\nolimits}

\newcommand{\Sym}{\mathop{\rm Sym}\nolimits}

\DeclareMathOperator{\GL}{\mathrm{GL}}
\DeclareMathOperator{\SL}{\mathrm{SL}}

\DeclareMathOperator{\diag}{\mathrm{diag}}

\DeclareMathOperator{\Span}{\mathrm{Span}}

\newcommand{\pr}{^{\,\prime}}

\theoremstyle{plain}
\newtheorem{thmS}[subsection]{Theorem}

\newtheorem{propS}[subsection]{Proposition}

\theoremstyle{definition}
\newtheorem{remS}[subsection]{Remark}
\newtheorem{defnS}[subsection]{Definition}

\theoremstyle{plain}
\newtheorem{thm}[subsubsection]{Theorem}
\newtheorem{lem}[subsubsection]{Lemma}
\newtheorem{prop}[subsubsection]{Proposition}
\newtheorem{cor}[subsubsection]{Corollary}

\theoremstyle{definition}
\newtheorem{rem}[subsubsection]{Remark}
\newtheorem{defn}[subsubsection]{Definition}

\newcommand{\sh}[1]{{#1}}

\def\ni{\noindent}

\begin{document}

\title{On Ideal generators for affine Schubert varieties}

\author{V. Kreiman}
\address{Department of Mathematics\\
Virginia Polytechnic Institute\\
Blacksburg, VA 24061}
\email{vkreiman@vt.edu}

\author[V. Lakshmibai]{V. Lakshmibai${}^{*}$}
\address{Department of Mathematics\\ Northeastern University\\
Boston, MA 02115}
\email{lakshmibai@neu.edu}
\thanks{${}^{*}$ Partially suported
by NSF grant DMS-0400679 and NSA-MDA904-03-1-0034.}

\author[P. Magyar]{P. Magyar${}^{\dag}$}
\address{Department of Mathematics\\
Michigan State University\\
East Lansing, MI 48824}
\email{magyar@math.msu.edu}
\thanks{${}^{\dag}$ Partially suported
by NSF grant DMS-0405948.}

\author[J. Weyman]{J. Weyman${}^{\ddag}$}
\address{Department of Mathematics\\
Northeastern University\\
Boston, MA 02115}
\email{j.weyman@neu.edu}
\thanks{${}^{\ddag}$ Partially
suported by NSF grant DMS-0300064.}

\maketitle

%\centerline{\Large\bf Standard Bases and ideal generation
%}\vspace{1em} \centerline{ V.~Kreiman} \centerline{ V.~Lakshmibai}
%\centerline{ P.~Magyar} \centerline{ J.~Weyman} \vspace{.5em}
%\centerline{June 2004} \vspace{0em}

\begin{abstract}
We consider a certain class of
Schubert varieties of the affine Grassmannian of type A.
By embedding a Schubert variety into a finite-dimensional
Grassmannian, we construct an explicit basis of sections
of the basic line bundle by restricting certain Pl\"ucker
co-ordinates.

As a consequence, we write an explicit set of generators for the
degree-one part of the ideal of the finite-dimensional embedding.
This in turn gives a set of generators for the degree-one part of
the ideal defining the affine Grassmannian inside the infinite
Grassmannian which we conjecture to be a complete set of ideal
generators.

We apply our results to the orbit closures of nilpotent matrices.
We describe (in a characteristic-free way) a filtration for the
coordinate ring of a nilpotent orbit closure and state a
conjecture on the SL(n)-module structures of the constituents of
this filtration.
\end{abstract}

\section{Introduction}
\setcounter{subsubsection}{0}

Let $K$ be the base field, which we shall
suppose to be algebraically closed of arbitrary characteristic.
Let $F:=K((t)),A:=K[[t]]$. Let
$\mathcal{G}=SL_n(F),\mathcal{P}=SL_n(A)$, so that
$\mathcal{G}/\mathcal{P}$ is the affine Grassmannian, and
let $X$ be an (affine) Schubert variety in
$\mathcal{G}/\mathcal{P}$.
In our previous paper \cite{klmw}, we constructed an explicit
basis for the affine Demazure module
$H^0(X,L)$, the space of sections of the basic line bundle $L$ on
$\mathcal{G}/\mathcal{P}$, restricted to $X$.
As a consequence, we obtained a basis for $H^0(X, L)$ in terms of
certain Pl\"ucker co-ordinates.

In this paper, a sequel to \cite{klmw},
we first give a different construction
of our basis of $H^0(X, L)$ for a certain class of affine Schubert
varieties.  Our key tool is the matrix
presentation of the elements of $\mathcal{G}/\mathcal{P}$,
similar to that for the classical Grassmannian.
This is constructed by embedding $\mathcal{G}/\mathcal{P}$
inside an infinite Grassmannian,
which in turn embeds each affine Schubert variety $X$
inside a certain finite-dimensional Grassmannian.

After defining the basis elements of $H^0(X, L)$
as restrictions of certain ``admissible" Pl\"ucker
co-ordinates (Definition \ref{adm}),
we prove that they form a spanning set using a system of
degree-one straightening relations, the
``shuffles" (Corollary \ref{shufflerels}).
We prove linear independence inductively
by writing certain degree-two relations
among the admissible Pl\"ucker coordinates
and restricting to a smaller $X$ (Section \ref{linind}).
Both the shuffles and the degree-two relations follow
from the matrix presentation.

As an important application, we obtain explicit generators for the
degree-one part of the ideal for our class of
affine Schubert varieties inside the finite-dimensional Grassmannian.
As a corollary, we obtain generators for
the degree-one part of the ideal of the
affine Grassmannian inside the infinite Grassmannian (Theorem
\ref{main2}); and we conjecture that these in fact give
a complete set of generators for the ideal. % (Section \ref{ele}).
Such ideal generators serve as effective tools
in solving the geometric problems, especially in
the study of the singularities.

It should be
remarked that though ideal generators for affine Schubert
varieties (as subvarieties of the affine Grassmannian)
may be found using Littelmann's standard monomial basis \cite{li},
they are hard to compute (since Littelmann's standard
monomial basis is hard to compute). Our goal is to develop a
standard monomial theory for affine Schubert varieties in terms of
Pl\"ucker co-ordinates, in the spirit of the classical work
of Hodge (\cite{ho1,ho2}; see also \cite{mu}). Construction of
a basis for the higher-level Demazure modules $H^0(X,L^{\otimes \ell})$
will be taken up in a subsequent paper.

Another application is to nilpotent orbit closures. Let $\NNN$
denote the set of all nilpotent matrices in $M_{n\times
n}(K)$, a closed affine subvariety of $M_{n\times n}(K)$.
The group $\GL_n(K)$ acts on $\NNN$ by conjugation. Each orbit
contains precisely one matrix in Jordan canonical form (up to
order of the Jordan blocks). Thus the orbits are indexed by
partitions $\mu$ of $n$, i.e. $\mu=(\mu_1,\ldots,\mu_n)$,
$n\geq\mu_1\geq\cdots\geq\mu_n\geq 0$, $\mu_1+\cdots+\mu_n=n$.
We denote the orbit corresponding to the partition
$\mu=(\mu_1,\ldots,\mu_n)$ by $\NNN_{\mu}^0$, and its closure by
$\NNN_{\mu}$. Lusztig has shown (cf.~\cite{lu}, \cite{lu2},
\cite{mag}) that each $\NNN_{\mu}$ is isomorphic to an
open subset of a certain affine Schubert variety
(generally not of the special class mentioned above).

Using Lusztig's isomorphism
(and the matrix representation for the elements of
$\mathcal{G}/\mathcal{P}$), we define in a characteristic-free
way a filtration ${\mathcal F}:=\{{\mathcal F}_{m,\mu},m\geq 0\}$ for $A_\mu$,
the co-ordinate ring of $\NNN_{\mu}$: namely,
${\mathcal F}_{m,\mu}$ is the span of monomials in the Pl\"ucker
co-ordinates on $\NNN_{\mu}$ of degree $\le m$.
\\[1em]
\textbf{Conjecture:} Let $E=K^n$. Let $\lambda = (\lambda_1
,\ldots ,\lambda_s )$ be the partition conjugate to $\mu$. Then
there is a characteristic free isomorphism of $SL(E)$-modules
\[{\mathcal F}_{m,\mu} \cong
L_{\lambda_1^m}E\otimes\ldots\otimes L_{\lambda_s^m}E\] where
$L_\nu E$ denotes the Weyl module with
highest weight $\nu$.
\\[1em]
In characteristic zero, this conjecture holds
due to \cite{mag2,Sh}.

The sections below are organized as follows: In \S \ref{infinite},
we describe the affine and the infinite Grassmannians, and the
embedding of $\mathcal{G}/\mathcal{P}$ inside the infinite
Grassmannian. In \S \ref{ind}, we describe the Ind-variety
structures for the affine and the infinite Grassmannians, and we
realize the affine Schubert varieties as closed subvarieties
inside certain finite-dimensional Grassmannians. In \S
\ref{sec_schubert _var}, we describe our special class of affine
Schubert varieties, and for $X$ in this class we establish our
main results concerning a basis for $H^0(X, L)$ and ideal
generators for $X$. We then apply these results to the affine
Grassmannian. In \S \ref{sec_nilpotent} we present the results for
the nilpotent orbit closures.

\section{Infinite and affine Grassmannians}
\label{infinite}
\setcounter{subsubsection}{0}

In this section, we recall the generalities on loop groups
and affine and infinite Grassmannians.
For details, we refer to \cite{kac,p-s} (see also \cite{mag}).
Let
$$
F:=K((t))=\left\{
\textstyle\sum\limits_{i=-N}^{\infty}a_it^i
\mid a_i\in K \right\}\,,
\quad
A:=K[[t]]=\left\{\textstyle\sum\limits_{i=0}^{\infty}a_i t^i\mid a_i\in K\right\}\,.
$$
If $f=\sum\limits_{i=-N}^{\infty}a_it^i\in F$, with $a_N\neq
0$, define $\ord(f):=N$, as well as $\ord(0):=-\infty$.
Let $n\in\NN$ be fixed.

\begin{defnS} $\Grhat(n)$, the  \emph{affine
Grassmannian}, is defined to be the set of all
$A$-lattices (free $A$-modules of rank $n$) in $F^n$.
\end{defnS}

\ni Define:
$$K^{\infty}:=\left\{
\textstyle\sum\limits_{j=N}^{\infty}a_je_j\mid \,N\in\ZZ,\,a_j\in K
\right\}\,,$$
and for $k\in\ZZ$, let:
$$E_k=\left\{\textstyle\sum\limits_{j=N}^{\infty} a_je_j\mid N\geq k,\ a_j\in K
\right\}\,.$$

\begin{defnS}  $\Gr(\infty)$, the  \emph{infinite
Grassmannian}, is defined to be the set of subspaces $V\subset
K^{\infty}$ such that $E_m\subset V\subset E_{-m}$, for some
$m\in\NN$.
\end{defnS}

Let $e_1,\ldots,e_n$ be the standard basis of $F^n$.
Consider the identification of $K$-vector spaces:
 $$\begin{array}{ccl}
 K^{\infty}&\simeq& F^n\\
 e_{cn+i}&\leftrightarrow& t^ce_i
 \end{array}\leqno(*)$$
 for $1\le i\le n$ and $c\in\ZZ$.
 We let $t:K^{\infty}\rightarrow K^{\infty}$ denote the map
 $e_i\mapsto e_{i+n}$,
 and we use the same symbol for the map
 $t:F^n\rightarrow F^n$ which is multiplication by $t$ .
 Indeed, these two maps are identified under ($*$).
 Also $\{V\in \Gr(\infty)\mid tV\subset V \}$
 gets identified with the set of $t$-stable
 subspaces of $F^n$, i.e., with the set of $A$-lattices in $F^n$ ;
 and we obtain an embedding $\Grhat(n)\subset\Gr(\infty)$.
 In particular, $E_1$ corresponds to the standard $A$-lattice
 $L_1=$ span$_A\{e_1,\cdots,e_{n}\}$.
% in particular, $E_i$ gets identified with the $A$-lattice
%$=\Span_A(e_i,e_{i+1},e_{i+2},\ldots,e_{i+{n-1}})$

For $V\in \Gr(\infty)$,
define the \emph{virtual dimension of} $V$ as
$$\vdim(V):=\dim(V\,/\,V\sh\cap E_1)-\dim(E_1\,/\,V\sh\cap E_1)\,.$$
For $j\in\ZZ$,
define $\Gr_j(\infty):=\{V\in \Gr(\infty)\mid \vdim(V)=j\}$\
and
$$
\Grhat_j(n):=\Grhat(n)\cap\Gr_j(\infty)
=\{V\in\Grhat(n)\mid\vdim(V)=j\}\,.
$$
We say that a subset $I\subset\ZZ$ is {\it almost natural} if
$|I\setminus I\cap\ZZ_+|$ and $|\ZZ_+\setminus I\cap\ZZ_+|$ are both
finite numbers, where $\ZZ_+=\{i\geq 1\}$.
In this case, we define the {\it virtual
cardinality of $I$} as:
$$\vcard(I)=\| I\|:=
|I\setminus I\cap\ZZ_+|-|\ZZ_+\setminus I\cap\ZZ_+|\,.$$
For such $I$, we have the associated coordinate subspace in
$\Gr(\infty)$, namely $$E_I=\left\{\sum_{i\in I}a_ie_i\mid a_i\in K\right\}\,,$$
and $\vdim(E_I)=\|I\|$.

If $(i_j)=(i_1<i_2<\cdots)$ is an increasing sequence of integers,
we say that $(i_j)$ is {\it almost natural} if the set $\{i_j\}$
is almost natural. In this case, we define
$\vcard\,(i_j):=\vcard\{i_j\}$. Observe that the increasing
sequence $(i_j)$ is almost natural with virtual cardinality $c$ if
and only if $i_j=j+c$ for $j\gg 0$. We denote the collection of
almost natural sequences by $\II$, and those of virtual
cardinality 0 by $\II_0$. We define the {\it Bruhat order} on
increasing sequences by:
$$
(i_j)\geq (k_j)\quad\Longleftrightarrow\quad
i_j\geq k_j\ \text{for all}\ j\,.
$$

Let $S_{\infty}$ be the group of permutations of $\ZZ$, and let
$\tau\in S_{\infty}$ be the permutation $\tau(i)=i+n$.  Define
the {\it extended affine Weyl group}
$\What:=\{w\in S_{\infty}\mid w\tau=\tau w\}$. Note that if
$w\in \What$, then $w$ is completely determined by $(w(1),\ldots,w(n))$.
Indeed, for $1\leq i\leq n$,
we have $w(nc+i)=w\tau^c(i)=\tau^c w(i)$.
Thus we can embed both the
finite permutation group $S_n$ and $\ZZ^n$ in $\What$ as follows.
We map $\sigma\in S_n$ to $w\in \What$ such that $w(i)=\sigma(i)$, $1\leq
i\leq n$; and we map $c=(c_1,\ldots,c_n)\in\ZZ^n$ to $w\in\What$ such
that $w(i)=\tau^{c_i}(r)$, $1\leq i\leq n$.  We identify $S_n$
and $\ZZ^n$ with their images under these embeddings.  Then
$\ZZ^n$ is a normal subgroup of $\What$, $\ZZ^n\cap S_n$ is the
identity, and $\mathbb{Z}^n S_n=\What$. That is, $\What=\ZZ^n\ltimes
S_n$, a semidirect product.

If $w\in \What$, then $\{w(j)\mid j\geq 1\}$ is almost natural.
Define $I_w$ to be the sequence obtained by listing the elements
of $\{w(j)\mid j\geq 1\}$ in increasing order. If we write $w=c\sigma$,
where  $c=(c_1,\ldots,c_n)\in\ZZ^n$ and $\sigma\in S_n$, then we may compute
 $\| I_w\|=-\sum c_i$.
%If $w_1$, $w_2\in \What$, we say that $w_1\geq w_2$ if $I_{w_1}\geq I_{w_2}$.
We use $\What^P$ to denote $\ZZ^n$ identified as
a set of coset representatives for $\What\!/S_n$.
% to be the maximal coset representatives in $\What$ of $\What/S_n$. Then
%$\What^P$ is precisely $\ZZ^n$.
One easily sees that a given $w\in \What^P$ is determined by its set
$I_w$, and that an almost-natural set $I$ can be realized as $I_w$
whenever $I$ is $\tau$-stable, meaining $\tau I\subset I$ (or equivalently
$tE_I\subset E_I$).  That is, we have a bijection
$$\begin{array}{rcl}
\What^P=\ZZ^n &\stackrel{\sim}{\rightarrow}&
\left\{I\subset\ZZ\left|\begin{array}{c}
 I\ \text{almost natural}\\ \tau I\subset I
 \end{array}\right.\right\}\\[1.5em]
w=(c_1,\ldots,c_n)&\mapsto& I_w = \{1\sh+c_1 n,\, 2\sh+c_2 n,\cdots,
n\sh+1\sh+c_1 n,\ldots\}\,.
\end{array}$$

Now define $W:=\{w\in\What\mid \|I_w\|=0\}$.
Then $S_n, \ZZ_0^n\subset W$, where $\ZZ_0^n:=\ZZ^n\cap
W=\{c\in\ZZ^n\mid c_1+\cdots+c_n=0\}$, and we have $W= \ZZ_0^n\ltimes S_n$,
a semidirect product.
Note that $W$ is the affine Weyl group associated to $\SLhat_n=\mathrm{SL}_n(F)$.
It is the Coxeter group generated by the adjacent transpositions
$s_1,\ldots,s_{n-1}\in S_n$ along with the reflection
$s_0$ defined by: $s_0(1)=0$,\ $s_0(n)=n\sh+1$, and $s_0(i)=i$ for $1<i<n$.
%Define $W^P$ to be the maximal coset representatives in $W$ of $W/S_n$.
%Then $W^P$ is precisely $\ZZ_0^n$.

We use $W^P$ to denote $\ZZ_0^n$ identified as
a set of coset representatives for $W/S_n$.
Once again $w\mapsto I_w$ gives a bijection:
$$
W^P=\ZZ_0^n\stackrel{\sim}{\to} \{I\in\II_0\mid
\tau I\subset I\}\,.
$$
We will denote $E_w:=E_{I_w}$.

Now, $\cg=\GL_n(K)$ acts
transitively on $\Grhat(n)$ and the isotropy subgroup at $E_1$ is
$\cp=\GL_n(A)$, so we have $\Grhat(n)\cong\cg/\cp$.
The group $\What$ embeds in $\cg$ : if $w=c\sigma$, where
$c=(c_1,\ldots,c_n)\in \ZZ^n$ and $\sigma\in S_n$, we identify $w$ with
the matrix $\diag(t^{c_1},\ldots,t^{c_n})\cdot [\sigma]\in \cg$,
where $[\sigma]$ is the permutation matrix associated to $\sigma$.
Let $\cb=\{g=(g_{ij})\in \GL_n(A)\mid
 \ord(g_{ij})>0\ \text{for}\ i<j\}$.
We have the Bruhat decomposition:
$$\cg=\prod\limits_{w\in \What^P}\cb w\cp\,,$$
and projection onto $\cg/\cp$ gives:
$$\Grhat(n)\cong\cg/\cp=\prod\limits_{w\in \What^P}\cb E_w,$$
where $E_w=w\cp\in \cg/\cp$, or equivalently $E_w = E_{I_w}\in\Grhat(n)$.
For any $V\in \cb E_w$, we have $\vdim(V)=\vdim(E_w)=\| I_w\|$. The
orbit $\cb E_w\subset \cg/\cp$ is called the {\it affine Schubert cell}
associated to $w$, and is denoted by $X^\circ(w)$.

Let $\cg_0:=\{g\in \cg\mid \ord(\det g)=0\}\supset\cp$.
Then $\cg_0$ acts transitively on $\Grhat_0(n)$ and the
isotropy subgroup at $E_1$ is $\cp$,
so $\Grhat_0(n)\cong\cg_0/\cp$.
We again have $W\hookrightarrow \cg_0$ and the Bruhat decomposition:
$$\cg_0=\prod\limits_{w\in W^P}\cb w\cp.$$
Projection onto $\cg_0/\cp$ gives:
$$\Grhat_0(n)\cong\cg_0/\cp=\prod\limits_{w\in W^P}\cb E_w\,.$$
Henceforth we focus on $\Grhat_0(n)$.
Note that $\SL_n(F)$ also acts transitively on
$\Grhat_0(n)$, the isotropy at $E_1$ being $\SL_n(A)$, so $\Grhat_0(n)\cong\SL_n(F)/\SL_n(A)$.

We will see in the following section that, although $\Grhat_0(n)\subset
\Gr_0(\infty)$ are Ind-varieties of infinite dimension, the Schubert
cell $X^{\circ}(w)$ is an ordinary variety with the finite dimension
$\sum_{j=1}^\infty (j-i_j)$, where $I_w = (i_1\sh<i_2\sh<\cdots)$
and $i_j=j$ for $j\gg0$.

\section{Ind-variety Structures}
\label{ind}
\setcounter{subsubsection}{0}

In this section, we recall the Ind-variety structure for affine
and infinite Grassmannians. See \cite{Ku,Ku2} for details.

For $s\in\NN$, let $V_s=t^{-s(n-1)}L_{1}/t^sL_1$, with
$\dim_K(V_s)=sn^2$. Consider the embedding
defined by:
$$\begin{array}{rcl}
\phi_s:\wedge^{sn}V_s&\hookrightarrow&\wedge^{(s+1)n}V_{s+1}\\[.5em]
v_1\wedge\cdots\wedge v_{sn}
&\mapsto&
v_1\wedge\cdots\wedge v_{sn}\wedge
t^s e_1\wedge\cdots\wedge t^s e_n
\,,
\end{array}$$
%t^{-(s+1)(n-1)}e_0\wedge\cdots\wedge t^{-s(n-1)-1}e_0\wedge
%v_1\wedge\cdots\wedge v_{sn}\wedge t^se_0
and let $\fuchs$ be the direct limit vector space $\underrightarrow{\lim} \,
(\wedge^{sn}V_{s})$.

The map  $\phi_s$ induces
$\phi_s^{\star}:\wedge^{(s+1)n}(V^{\star})_{s+1}\to
\wedge^{sn}(V^{\star})_{s}$, and we define $\fuchs^{\star}$ to be the
vector space $\underleftarrow{\lim} \, (\wedge^{sn}V^{\star}_s)$.
There is a bilinear pairing between $\fuchs$ and $\fuchs^{\star}$
implied by the universal properties of limits under which
$\fuchs^{\star}=(\fuchs)^{\star}$, the dual space of $\fuchs$.
The universal properties of limits also imply
for all $s$ the existence of an injection $i_s:\wedge^{sn}V_s\to
\fuchs$ and a projection $\pi_s:\fuchs^{\star}\to \wedge^{sn}
(V^{\star}_s)=(\wedge^{sn}V_s)^{\star}$.

If $I=(i_1,i_2,\ldots)\in\II_0$, we can view $e_I:=e_{i_1}\wedge
e_{i_2}\wedge\cdots$ as an element of $\fuchs$ and
$p_I:=e^{\star}_{i_1}\wedge e^{\star}_{i_2}\wedge\cdots$ as an
element of $\fuchs^{\star}$. We call the $p_I$ the
{\it infinite Pl\"ucker coordinates}.  Indeed $\{e_I\mid
I\in\II_0\}$ forms a basis for $\fuchs$ and $\{p_I\mid I
\in\II_0\}$ is the dual basis for $\fuchs^{\star}$, i.e., $\langle
e_I,p_J\rangle=\delta_{IJ}$ (Kronecker delta) for
$I,J\in\II_0$.

The map $\phi_s$ descends to a closed immersion
$\PP(\wedge^{sn}V_s)\to\PP(\wedge^{(s+1)n}V_{s+1})$~, and we let
$\PP(\fuchs)=\bigcup_{s\in\NN} \PP(\wedge^{sn}V_s)$ with the
reduced projective Ind-variety structure (see \cite{Ku},
\cite{Ku2}). The Ind-variety homogeneous coordinate ring of
$\PP(\fuchs)$ is defined to be
$K[\PP(\fuchs)]:=\underleftarrow{\lim} \,
K[\wedge^{sn}V_s]=\Sym(\fuchs^{\star})$.

Fix an integer $s\geq 0$, and let
$$
\FF_s=\left\{V\in\Gr_0(\infty)
\left|\begin{array}{c}
t^{-s(n+1)}E_{1}\supset V\supset t^sE_{1} \\
\dim_K(V/t^sE_1)=sn
\end{array}\right.\right\}\,.
$$
%where $S=(-(r-s)n,-(r-s-1)n,\cdots,0,
%n,2n,\cdots,sn(=r),r+1,r+2,\cdots)$.
%Note that under ($*$), $E_{sn}$
%gets identified with the $A$-lattice $t^sL_1$ and $E_S$ with
%$t^{-s(n-1)}L_0$.
We identify $\FF_s$ with its image under the
bijection $b_s:\FF_s\to\Gr(sn,V_s)$,\ $V\mapsto V/t^sE_1$.
Now $\Gr(sn,V_s)$ embeds as a closed subvariety
of $\PP(\wedge^{sn}V_s)$ under the Pl\"ucker embedding $j_s$, and we
have a commutative diagram:
\begin{displaymath}
\begin{CD}
\Gr(sn,V_s)   @>j_s>>     \PP(\wedge^{sn}V_s)\\
@VVV                               @VV\phi_s V \\
\Gr((s+1)n,V_{s+1})           @>j_{s+1}>>
\PP(\wedge^{(s+1)n}V_{s+1})\,.
\end{CD}
\end{displaymath}
Thus
$\bigcup\limits_{s\geq 0}\FF_s\cong\bigcup\nolimits_{s\geq 0}\Gr(sn,V_s)$
induces the structure of a closed Ind-subvariety of
$\PP(\fuchs)$.

We claim that $\Gr_0(\infty)=\bigcup\nolimits_{s\geq 0}\FF_s$.
Indeed, for $V\in\Gr(\infty)$ with
$t^{-s(n+1)}E_{1}\supset V\supset t^sE_{1}$,
we have $V\in\FF_s$ if and only if
\begin{eqnarray*}
sn&=&\dim(V/t^sE_1)\\
&=&\dim(V\,/\,V\sh\cap E_1)+\dim(V\sh\cap E_1\,/\,t^sE_1)\\
&=&\dim(V\,/\,V\sh\cap E_1)+(\dim(E_1/t^sE_1)-\dim(E_1\,/\,V\sh\cap E_1))\\
&=&\dim(V\,/\,V\sh\cap E_1)+sn-\dim(E_1\,/\,V\sh\cap E_1),
\end{eqnarray*}
i.e., whenever $\vdim(V)=0$, meaning $V\in\Gr_0(\infty)$.

Multiplication by $t$ induces
%an endomorphism on $K^{\infty}$, which induces
a nilpotent endomorphism $t_s$ on each $V_s$. Define
$u_s:=1+t_s$, a unipotent automorphism of $V_s$, and denote
the induced automorphism of $\Gr(sn,V_s)$ also by $u_s$.
Let $$
\HH_s=\left\{V\in \Grhat_0(n)\left|
\begin{array}{c}
 t^{-s(n-1)}L_{1}\supset V\supset t^sL_1\\
 \dim_K(V/t^sL_1)=sn\end{array}\right.\right\}\,.
 $$
 Using the identification $(*)$ of Section 2 and
 the isomorphism $b_s$ above,  we have:
 $$
 \HH_s\ \stackrel{(*)}\cong\
 \{V\in\FF_s\mid tV\subset V\}
 \ \stackrel{b_s}\cong\
 \Gr(sn,V_s)^{u_s}\,,$$
a closed projective subvariety of $\Gr(sn,V_s)$.
Therefore
$\Grhat_0(n)=\bigcup\nolimits_{s\geq 0}\HH_s$ $\cong\bigcup\nolimits_{s\geq 0}
\Gr(sn,V_s)^{u_s}$ induces the structure of a closed
Ind-subvariety of $\Gr_0(\infty)$.
For $w\in W^P$, we define the {\it affine Schubert variety}
$X(w)$ to be the Zariski closure
$\overline{X^{\circ}(w)}=\overline{\cb E_w}\subset\Grhat_0(n)$.

The following diagram illustrates the relationships between the
various varieties and Ind-varieties:
\begin{equation}\label{diagprojections1}
\begin{CD}
X(w) @>>> \Grhat_0(n) @>>> \Gr_0(\infty)  @>>> \PP(\fuchs)\\
@. @AAA               @AAA                 @AAA   \\
@. \Gr(sn,V_s)^{u_s} @>>> \Gr(sn,V_s)  @>>>  \PP(\wedge^{sn}V_s)
\end{CD}
\end{equation}
If we consider each projective variety in the bottom row as a
closed Ind-subvariety of the projective Ind-variety above it, all
maps in the diagram are closed immersions of Ind-varieties.

Thus we have the corresponding graded projections (i.e. restrictions)
of the homogeneous coordinate rings:
\begin{equation}\label{diagprojections2}
\hspace{-1em}
\begin{CD}
@. @. @. \Sym(\fuchs^{\star})\\
@. @. @. ||\\
K[X(w)] @<<< K[\Grhat_0(n)] @<<< K[\Gr_0(\infty)]  @<<< K[\PP(\fuchs)]\\
@. @VVV               @VVV                 @VVV   \\
@. K[\Gr(sn,V_s)^{u_s}] @<<< K[\Gr(sn,V_s)]  @<<<
K[\PP(\wedge^{sn}V_s)]\\
@. @. @. ||\\
@. @. @. \Sym((\wedge^{sn}V_s)^{\star})\\
\end{CD}
\end{equation}
%\begin{equation}\label{diagprojections2}
%\begin{CD}
%K[\Grhat_0(n)] @<<< K[\Gr_0(\infty)]  @<<< K[\PP(\fuchs)]=\Sym(\fuchs^{\star})\\
%@VVV               @VVV                 @VVV   \\
%K[\Gr(n,n^2)^{u_n}] @<<< K[\Gr(n,n^2)]  @<<<
%K[\PP(\wedge^nV_n)]=\Sym((\wedge^nV_n)^{\star})
%\end{CD}
%\end{equation}
Note that the projections of the infinite Pl\"ucker coordinates to
the third row of (\ref{diagprojections2}) are just the usual
finite Pl\"ucker coordinates.

\begin{propS}\label{pluckervanish}
\begin{itemize}
\item[(i)] $p_I|_{X(w)}=0\iff I\not\leq I_w$
%\item $X(w)$ is the vanishing set in $\Grhat_0(n)$ of the $P_I$,
%$I\in\II_0$, $I\not\leq I_w$.
\item[(ii)] For $w\in W^P$, we have $X(w)=\dot{\cup} X^{\circ}(y)$,
the disjoint union
being over $\{y\in W^P\mid y\leq w\}$.
\item[(iii)] For $
w_s=(c_1,\ldots,c_n):=(-s(n-1),s,\ldots,s)\in W^P\,,
$
we have the isomorphisms of reduced varieties:
$$
X(w_s)=\left\{
\text{\rm $A$-lattices}\ L\,
\left|\,\begin{array}{c}
t^{-s(n-1)}L_1\supset L\supset t^sL_1\\
\dim(L/t^sL_1)=sn
\end{array}\right.\right\}
\ \cong\ \Gr(sn,V_s)^{u_s}\,.
$$
\end{itemize}
\end{propS}
\begin{remS}\label{cano} The Proposition implies:\ \
$\Grhat_0(n)=\underrightarrow{\lim} \, X(w_s)$.
\end{remS}

We will adopt an abuse of notation when dealing with
almost natural sets.
For any $I=(i_1,i_2,\cdots)\in \II_0$,
there is some $m\geq 0$ such that $i_j=j$ for
$j\geq m$.   We shall ignore the ``trivial'' entries
$i_j=j$ for $j\geq m$, and write
just the $(m\sh-1)$-tuple $I=(i_1,\cdots,i_{m-1})$.
Any such sequence may also be considered as an $\ell$-tuple
for $\ell\geq m\sh-1$ by restoring some of the entries $i_j=j$.
Thus, any $I\in\II_0$ may be thought of as an $sn$-tuple for $s\gg 0$.
Now part (ii) of the Proposition implies that
every $X(w)$ sits inside $X(w_s)$ for $s\gg 0$.

The action by left translations of $\SL_n(F)$ on $\Grhat_0(n)$
induces a projective representation of $\SL_n(F)$ and its
Lie algebra $\mathfrak{sl}_n(F)$ on $K[\Grhat_0(n)]$.
This lifts to an actual representation of a central extension,
the affine Kac-Moody algebra $\widehat{\mathfrak{sl}}_n$
acting on $K[\Grhat_0(n)]$.

\begin{thmS} (cf. \cite{Ku2})\label{th_irrrep_coordring}
\begin{itemize}
\item[(i)] Let $V_{\Lam_0}$ be the basic integrable irreducible
representation of $\widehat{\mathfrak{sl}}_n$.  Then
$V^{\star}_{d\Lam_0}\cong
K[\Grhat_0(n)]_d=\Sym^d({\fuchs^{\star}})|_{\Grhat_0(n)}$ as
$\widehat{\mathfrak{sl}}_n$ representations.
\item[(ii)] For $w\in W^P$ let $V_{\Lam_0}(w)$ be the affine
Demazure module of $V_{\Lam_0}$ corresponding to $w$. Then
$$V^{\star}_{d\Lam_0}(w)\cong
K[X(w)]_d=\Sym^d({\fuchs^{\star}})|_{X(w)}\,.$$
\end{itemize}
\end{thmS}

\section{The Schubert variety $X(w_s)$}\label{sec_schubert
_var}
\setcounter{subsubsection}{0}

\subsection{Renormalization of indices}

Recall the special element $$w_s=(-s(n-1),s,\ldots,s)\in W^P=\ZZ^n_0\,,$$
with
$
I_{w_s}=(1-s(n\sh-1)n,\ 1-s(n\sh-1)n+n,\ldots,1+sn-n,\ 1+sn)\,.
$

In order to avoid continually writing negative indices, we
now renormalize these values.  This means that,
instead of working with the component $\Grhat_0(n)$, we consider
the isomorphic Ind-variety
$$
\Grhat_0\pr(n):=t^{s(n-1)}\Grhat_0(n)\ =\
\{ V\in\Grhat(n) \mid \vdim(V)=-s(n-1)n\}\,.
$$
The Bruhat decomposition becomes:
$$
\Grhat_0\pr(n) = \bigcup_{w} X^\circ(w)\,,
$$
where the union is over $w$ in:
$$\begin{array}{rcl}
\Acute W^P&:=&\tau^{s(n-1)}W^P = \tau^{s(n-1)}\ZZ^n_0\\[.5em]
&=& \{\ w=(c_1,\ldots,c_n)
\ | \ c_1\sh+\cdots\sh+c_n= s(n\sh-1)n\}\,,
\end{array}$$
so that
$$I_w\in \II_0\pr\ :=\ \tau^{s(n-1)}\II_0\ =\ \II_{-s(n-1)n}\,.$$

In this normalization,
the zero-dimensional Schubert variety is $X(w)=\{E_w\}$
for $$
w=\tau^{s(n-1)}
=(s(n-1),\ldots,s(n-1))
$$
$$
I_{w}=\{1\sh+sn(n\sh-1),2\sh+sn(n\sh-1),
3\sh+sn(n\sh-1),\ldots\}\,.
$$
We renormalize:
$$
w_s:=(0,sn,\ldots,sn)
$$
$$
I_{w_s}=(1,n\sh+1,2n\sh+1,\cdots
,sn^2\sh+1,sn^2\sh+2,sn^2\sh+3,\ldots)\,.
$$
Further, we can rephrase Proposition \ref{pluckervanish}(iii) as:
$$
X(w_s)=
\left\{\text{$A$-lattices}\ L\subset F^n\
\left|\begin{array}{c}
L_1\supset L\supset t^{sn}L_1\\[.5em]
\dim L/t^{sn}L_1=sn\end{array}\right.\right\}
\cong
\Gr(sn,V_s)^{u_s}\,,
$$
where $V_s= L_1/t^{sn}L_1$.

\subsection{$\ZZ\sh\times\ZZ_+$ matrix presentation} \label{ele}

As with the finite-dimensional Grassmannians, an element $V \in
\Gr(\infty)$ may be thought of as the column space of a matrix
with two-sided infinite columns and one-sided infinite rows: a
$\ZZ\sh\times\ZZ_+$ matrix $M=(a_{ij})_{(i,j)\in\ZZ\times\ZZ_+}$.
Since $E_j\subset V$ for $j\gg 0$, we write $A$ so that
$a_{ij}=\delta_{ij}$ except for finitely many entries. Using the
$(*)$ identification $K^{\infty}\simeq F^n,\
e_{cn+i}\leftrightarrow t^ce_i$, an $A$-lattice in $F^n$ also has
a $\ZZ\sh\times\ZZ_+$ matrix presentation.

For example, we have the coordinate lattice:
$E_I=\Span_K(e_i\mid i\in I)$ for $I\in \II$.
Writing each basis element $e_i$
as a column vector with a 1 in row $i$ and 0 elsewhere,
we obtain the matrix presentation of $E_I$.
%For $w=(c_1,\ldots,c_n)\in \What^P=\ZZ^n$,
%recall that $I_w$ is the increasing rearrangement
%of $\{1\sh+c_1n,\ 2\sh+c_2n,\ldots, n\sh+1\sh+c_1n,\ n\sh+2\sh+c_2n,\ldots\}$,
%and $E_w=E_{I_w}$.

%The open cell $X^\circ(w)$ consists
%of  $\ZZ\sh\times\ZZ_+$ matrices with (the relevant) columns being
%determined by the entries of $S_w$; further, in a typical column,
%the entries above the diagonal entries are zero (by our choice of
%the Borel sub group).
%In this section, we are interested in
Now let us consider $w=w_s$,\ $I=I_{w_s}$.  We set:
$$
\fbox{$\begin{array}{c}r:=sn\end{array}$}\ ,
$$
so that $X(w_s)\subset\Gr(r,V_s)$ where
$$\begin{array}{rcl}
V_s=\Span_K\{t^{c}e_1,\ldots,t^c e_{n}\mid 0\leq c\leq r\sh-1\}\,.
\end{array}$$
We may write each vector $v\in V_s$ with respect to
this basis in the block form:
\begin{displaymath}
v=\left(\begin{array}{c} v_1\\v_2\\ \vdots \\ v_{r}
\end{array}\right)
\end{displaymath}
where each $v_i$ is an $n\sh\times 1$ column
vector.  Multiplication by $t$ becomes:
\begin{displaymath}
t(v)=\left(\begin{array}{c}
0\\v_1\\v_2\\ \vdots \\ v_{r-1}
\end{array}\right)
\end{displaymath}

Taking the $rn\times r$ matrix presentation of $\Gr(r,V_s)$,
we identify
$$
\Gr(r,V_s)\
\cong\ M_{rn\times r}^{\max}(K)\,/\, \GL_{r}(K)\,,
$$
where $M^{\max{}}$ indicates matrices of maximal rank,
and two matrices define the same subspace $V\subset V_s$
if they differ by a column transformation in $\GL_{r}(K)$.
The $rn\times r$ matrix $M$ of a space $V\in \Gr(r,V_s)$
is the ``relevant part" of the $\ZZ\times\ZZ_+$ matrix corresponding to
the lattice $L = V\oplus t^rL_1\in \Grhat(n)$.

For a generic subspace $V$ in the Schubert cell $X^{\circ}(w_s)
\subset\Gr(r,V_s)^{u_s}$, we can find a {\it generating vector} $v\in V$
such that its $e_1$ component is equal to 1 and
$V=\Span_{k}(v,tv,\ldots,t^{r-1}v)$:  indeed, any vector
$v=e_1+a_2v_2+\cdots$ is a generating vector.
Thus $V$ can be presented by an $rn\times r$ matrix of the
form:
$$
M=(a_{ij})=\left(\begin{array}{cccc} v & tv & \cdots & t^{r-1}v
\end{array}\right)\ =\
\left(\begin{array}{c}
A_1\\ A_2\\ \vdots \\ A_r
\end{array}\right)\,,$$
where each $A_k$ is an $n\times r$ block.
This is a lower triangular $rn\times r$ matrix,
meaning $a_{ij}=0$ for $i<j$ ,
and the $1^{{\text{st}}}$ column of $A_1$
repeats as the $2^{{\text{nd}}}$
column of  $A_2$, the $3^{{\text{rd}}}$
column of  $A_3$, etc.

In particular, we have
 $1=a_{1,1}=a_{1+n,2}=\cdots=a_{(r-1)n+1,r}$.
 % $a_{(m-1)n+1,m}=1$ for $m=1,\ldots,r$.
We shall refer to these rows $1,\ 1\sh+n,\ldots,\,1\sh+(r\sh-1)n$
 % $((m-1)n+1)^{{\text{th}}}$ for  $m=1,\ldots,r$
as the {\it pivotal rows}.
Now we may reduce $M$ by column
operations so that each $a_{1,1}=a_{n+1,2}=\cdots=1$
%$a_{(m-1)n+1,m}$
is the {\it only} nonzero entry in a pivotal row.
This normalizes $M$ to make it a unique represenative
of $V\in X^\circ(w_s)$.
\vfill\eject
\ni That is, we may identify $X^{\circ}(w_s)\subset\Gr(r,V_s)$ with
the affine space of matrices $M=(a_{i,j})\in M_{rn\times r}(K)$ of the form:
\begin{displaymath} M=\left[\begin{array}{cccccc}
1&0&0&0&\cdots&0\\
a_{2,1}&0&0&0&\cdots&0\\
a_{3,1}&0&0&0&\cdots&0\\
\vdots&\vdots&\vdots&\vdots&\vdots&\vdots\\
a_{n,1}&0&0&0&\cdots&0\\
%\\
\hline%\\
0&1&0&0&\cdots&0\\
a_{n+2,1}&a_{2,1}&0&0&\cdots&0\\
a_{n+3,1}&a_{3,1}&0&0&\cdots&0\\
\vdots&\vdots&\vdots&\vdots&\vdots&\vdots\\
a_{2n,1}&a_{n,1}&0&0&\cdots&0\\
%\\
\hline%\\
0&0&1&0&\cdots&0\\
a_{2n+2,1}&a_{n+2,1}&a_{2,1}&0&\cdots&0\\
a_{2n+3,1}&a_{n+3,1}&a_{3,1}&0&\cdots&0\\
\vdots&\vdots&\vdots&\vdots&\vdots&\vdots\\
a_{3n,1}&a_{2n,1}&a_{n,1}&0&\cdots&0\\
%\\
\hline%\\
%\vdots&\vdots&\vdots&\vdots&\vdots&\vdots\\
\vdots&\vdots&\vdots&\vdots&\vdots&\vdots\\
%\\
\hline%\\
0&0&0&0&\cdots&1\\
a_{(r-1)n+2,1}&a_{(r-2)n+2,1}&a_{(r-3)n+2,1}&\cdots&\cdots&a_{2,1}\\
a_{(r-1)n+3,1}&a_{(r-2)n+3,1}&a_{(r-3)n+3,1}&\cdots&\cdots&a_{3,1}\\
\vdots&\vdots&\vdots&\vdots&\vdots&\vdots\\
a_{rn,1}&a_{(r-1)n,1}&a_{(r-2)n,1}&\cdots&\cdots&a_{n,1}\\
\end{array}\ .
\right]\end{displaymath}
%(note: indices
%$1,{n+1},\ldots,{(s-1)n+1}$ have been omitted for notational
%convenience). The horizontal lines in the above matrix divide $M$
%into $sn$ sub-matrices of size $n\times sn$ which we call {\it
%blocks}.

Now it is clear from inspecting the above matrix that,
the infinite Pl\"ucker coordinates on $\Gr_0(\infty)$,
when restricted to $\Gr(r,V_s)$, either vanish or
become usual Pl\"ucker coordinates on $\Gr(r,V_s)$,
the $r\times r$ minors of the above matrix $M$.
Indeed, the coordinate $p_I|_{X(w_{s})}$ is non-zero
for $I=(i_1<i_2<\cdots)$ if and only if
\begin{enumerate}
\item $i_j\ge (j\sh-1)n\sh+1$\ \ for\ \ $j=1,2,\ldots,r$.
\item $i_j=j+rn$\ \ for\ \ $j\geq r+1$.
\end{enumerate}
We shall denote such an index set as an $r$-tuple
$I=(i_1,\ldots,i_{r})$.

%\begin{rem}\label{r-tuple} In the sequel, we shall identify $w_s$ with the
%$r$-tuple $(1,n+1,2n+1,\cdots ,(r-2)n+1,(r-1)n+1)$. With this
%convention, we have that for a $S=(s_1,s_2,\cdots ),
%p_{S}|_{X(w_{s})}$ is non-zero if and only if
%\ni 1. $s_j\ge (j-1)n+1,1\le j\le r$
%\ni 2. $\{s_{r+1},s_{r+2},\cdots \}$ is the natural sequence
%$\{rn+1,rn+2,\cdots\}$. (Note that for the above normalization, we
%have that if $M$ represents a generic element in $X^{\circ}(w_s)$, then
%in the $(rn+j)$-th row of $M,j=1,2,\cdots $, the $(r+j)$-th entry
%is $1$ and all other entries are $0$.)
%\ni We shall denote such a $S$ by just $(s_1, \cdots ,s_r)$; then
%the corresponding $p_S$'s are precisely the Pl\"ucker coordinates
%on $\Gr(sn,V_s)$ which do not vanish on $X(w_s)$. In particular, a
%$\tau\in W^P,\tau\le w_s$ will be denoted as a $r$-tuple (note
%that Id corresponds to the (largest)$r$-tuple
%$((r-1)n+1,(r-1)n+2,\cdots,rn)$).
%\end{rem}

\subsection{Generation by admissible Pl\"ucker Coordinates}

\begin{defn}\label{adm}
A set $I=(i_1,i_2,\ldots)\in\II_0\pr$ is {\it admissible} if
$$
i_{j+1}-i_j\leq n\ \ \text{for}\ \ j\geq 1\,.
$$
In this case we also say that the Pl\"ucker coordinate
$p_I\in K[\Grhat_0\pr(n)]$ is admissible.
If $w\in \Acute W^P$, then $I$ (or $p_I$) is said to be {\it admissible} on
$X(w)$ if $I$ is admissible and $p_I|_{X(w)}\neq 0$.
\end{defn}
% equivalently, $P_I$
%(or $I$) is admissible on $X(w)$ if $P_I$ is admissible and $I\leq
%I_w$ (see Proposition \ref{pluckervanish} (i)).  Let $w$, $w'\in
%W^P$, and let $w'\leq w$. By definition, the admissible Pl\"ucker
%coordinates on $X_{w'}$ are the admissible Pl\"ucker coordinates
%on $X(w)$ which do not vanish on $X_{w'}$.
%\begin{rem}\label{t-stability}
%For $w\in W^P$, we have, $I_w$ has the property that if some entry
%$a$ is in $I_w$, then $a+jn$ is in $I_w$, for all $j=1,2,\ldots$;
%in fact for $I\in \II_0$, we have, $I=I_w$ for some $w\in W^P$ if
%and only if $I$ has the property that if some entry $a$ is in $I$,
%then $a+jn$ is in $I$, for all $j=1,2,\ldots$ - we shall refer to
%this as {\it $t$-stability} for $I$. We have, if $I=I_w$, then
%$I_w$ is admissible necessarily.
%\end{rem}

We proceed to show that on $X(w_s)$,
a non-admissible Pl\"ucker coordinate
$p_I$ is a linear combination of admissible Pl\"ucker coordinates.

\begin{lem}\label{rowcoltrick}
Let $A, B\in M_{m\times m}(K)$. Fix $\ell\leq m$.  Then the
sum of all determinants obtained by replacing $\ell$ rows of $A$ by
the corresponding $\ell$ rows of $B$ equals the sum of all
determinants obtained by replacing $\ell$ columns of $A$ by the
corresponding $\ell$ columns of $B$.
\end{lem}
\begin{proof} One sees by performing Laplace expansions that (1) the
sum of all determinants obtained by replacing $\ell$ rows of $A$ by
the corresponding $\ell$ rows of $B$ ; and (2) the sum of all determinants
obtained by replacing $\ell$ columns of $A$ by the corresponding $\ell$
columns of $B$ ; are both equal to (3) the sum of all products of $\ell$
minors of $B$ with complementary $(m-\ell)$ minors of
$A$.
\end{proof}

Let $I=(i_1,\ldots,i_{r})$ be such that
$p_I|_{X(w_{s})}\not=0$, so that $i_{r}\le
rn$.  Further let $i_1\geq n\sh+1$,
and let $\widetilde I=(i_1\sh-n,\ldots,i_{r}\sh-n)$. Choose $\ell\leq r$.
We define a {\it shuffle} to be
the sum of all Pl\"ucker coordinates obtained by
replacing $\ell$ elements $J\subset I$ with the corresponding $\ell$
elements $\widetilde J\subset\widetilde I$:
$$
\mathrm{sh}_I:=\mathop{\sum_{J\,\subset\, I}}_{|J|=\ell}
p_{I\,\setminus\, J\,\cup\, \widetilde J} \,.
$$

\begin{cor}\label{shufflerels}
The shuffles are identically $0$ as functions on $X(w_s)$:
$$\mathrm{sh}_I|_{X(w_s)} = 0\,.$$
\end{cor}

\begin{proof}  Let $I,\tilde I$ be as above.  Let $M_{rn\times r}$ be
the normalized matrix corresponding to a generic point of
$X^\circ(w_s)$, and let $M$ (resp. $\widetilde M$) be the $r\times
r$ submatrix of $M_{rn\times r}$ with row indices given by $I$
(resp. $\widetilde I$). Then we have that the last $r\sh-1$
columns of $M$ are same as the first $r\sh-1$ columns of
$\widetilde M$, and the last column of $\widetilde M$ consists of
zeroes. (Note that since $i_{r}\le rn$, we have, $i_{r}\sh-n\le
(r\sh-1)n$. Hence, the entries of $\widetilde M$ are taken from
the first $(r-1)$ blocks of $M_{rn\times r}$.)

This implies the vanishing of the minor obtained by replacing any $\ell$ columns
of $M$ with the corresponding $\ell$ columns of $\widetilde M$. Hence
in view of Lemma \ref{rowcoltrick}, we obtain the vanishing of the sum of all
minors of $M$ obtained by replacing $\ell$ rows of $M$ with the
corresponding $\ell$ rows of $\widetilde M$. But this sum is simply the
evaluation on $M_{rn\times r}$ of the sum of all Pl\"ucker coordinates obtained
by replacing $\ell$ indices of $I$ with the corresponding $\ell$
indices of ${\widetilde I}$. The desired result now follows.
\end{proof}

\begin{prop}\label{gen}
Let $I=(i_1, \cdots, i_r)$ be a non-admissible $r$-tuple such that
$p_I\,|\,_{X(w_{s})}\not=0$. Then on $X(w_s)$ we have,
$p_I=\sum\nolimits_{S}a_S p_S$ ,\ $a_S\in K$, where the sum is over
admissible $r$-tuples.
\end{prop}

\begin{proof}
Let $\ell\le r\sh-1$ be the smallest integer such that
$i_{\ell+1}-i_\ell>n$. Let $$I_1=(i_1+n,\cdots,i_\ell+n,
i_{\ell+1},\cdots,i_r),\ I_2=(i_1,\cdots,i_\ell,
i_{\ell+1}-n,\cdots,i_r-n)$$
Note that the (positive) difference
between the respective entries of $I_1,I_2$ equals $n$. Hence by
Corollary \ref{shufflerels}, the corresponding shuffle equals $0$
on $X(w_s)$. Let $p_J$ be a typical term appearing in the shuffle.
Corresponding to the switch of $\{i_1\sh+n,\cdots,i_\ell\sh+n\}$
respectively with $\{i_1,\cdots,i_\ell\}$, we have, $J$ equals $I$,
and any other $J$ is lexicographically greater than $I$. Now the
result follows by induction: note that the lexicographically
largest $r$-tuple is $((n\sh-1)r+1,\,(n\sh-1)r+2,\cdots,rn)$ which
corresponds to the identity or zero-dimensional Schubert variety).
\end{proof}
%Next we introduce new notation in order to reformulate the shuffle
%relations in a manner which will be useful in Chapter
%\ref{ch_thin}.
%Let $I=i_1, \ldots, i_n$ be an n-tuple, and let $H, J\subset I$.
%Define $\HI= \Hi_1, \ldots, \Hi_n$ and $\JI= \Ji_1, \ldots, \Ji_n$
%as follows:
%\begin{displaymath}
%\Hi_m=\left\{\begin{array}{ll} i_m & m\not\in H\\ i_m+n & m\in H
%\end{array}\right.
%\end{displaymath}
%\begin{displaymath}
%\Ji_m=\left\{\begin{array}{ll} i_m & m\not\in J\\ i_m-n & m\in J
%\end{array}\right.
%\end{displaymath}
%$m=1,\ldots, n$.  Define $\JHI:= \, _J(\HI)$.
%\begin{prop}\label{linstraight}
%Let $H\subset I$.  Then $\DS\sum\limits_{\begin{array}{l} J\subset I\\
%|J|=|H|
%\end{array}}P_{\JHI}=0$ on $X(w_s)$.
%\end{prop}
%\begin{proof} Consider the n-tuples $\HI=\Hi_1, \ldots, \Hi_n$
%and $\HI-n=\Hi_1-n, \ldots, \Hi_n-n$.  Then $\{ \JHI\mid J\subset
%I, |J|=|H|\}$ is the set of all n-tuples obtained from $\HI$ by
%replacing $|H|$ elements of $\HI$ by the corresponding $|H|$
%elements of $\HI-n$. The rest of the proof follows in a similar
%fashion as the proof of Corrolary \ref{shufflerels}.\end{proof}

\subsection{The reduced chain of $w_s$}

As an element of the Weyl group $W=\langle s_0,s_1,\ldots,s_{n-1}\rangle$,
we can write a reduced expression for $w_s$ as:
$$w_s={\underbrace{(s_1s_2\cdots s_{n-1}s_0)\cdots
(s_1s_2\cdots s_{n-1}s_0)}_{s(n-1)\ {\text{times}}}}$$ (In
particular, $\dim_K X(w_s)=s(n\sh-1)n=r(n\sh-1)$, which is indeed
the number of variables in the generic matrix above.)

Consider the full chain of $w$'s obtained by deleting one simple
reflection at a time in the above reduced expression for $w_s$
(starting from the left). Let us denote these by $\tau_i$ for
$i=1,\ldots, s(n-1)n+1=r(n-1)+1,$ so that
$$w_s=\tau_1>\tau_2>\cdots>\tau_{r(n-1)+1}=\mathrm{id}\,.$$
% dim $X(\tau_i)=s(n-1)n-i$; in
%particular we have, $\tau=w_s,\ \tau_{s(n-1)n}=id
%(=\tau_{r(n-1)+1})$.
Let us collect these into $r$ groups of $n$ elements so that
\ni the first group consists of $\{\tau_1,\cdots ,\tau_n\}$,
\ni the second group consists of $\{\tau_n,\cdots ,\tau_{2n-1}\}$,
\ni the next group consists of $\{\tau_{2n-1},\cdots
,\tau_{3n-2}\}$, etc.: that is, each group has $n$ elements and the first element
in each group is the last element of the previous group.
%(Note that this accounts for $rn-(r-1)=r(n-1)+1$
%(= the number of elements in our chain.)

The tuples corresponding to the $n$ elements $\{\tau_1,\cdots
,\tau_n\}$ in the first group are given by:
$$\tau_j=(j,n\sh+j,\ldots ,(r\sh-1)n\sh+j)\ \ \text{for}\ \ 1\le j\le n\,.$$
In particular, we have:
$$\tau_n=(n,2n,\cdots ,rn)\,.$$

More generally, for $2\le i\le r$,
let us denote the $i$-th group (of $n$ elements) by
$$\varphi_{ij}:=\tau_{(i-1)n-(i-j)+1}\,,\quad 1\le j\le n\,.$$
For $k\sh<l$, let $[k,l]:=(k,k\sh+1,\cdots ,l)$.
Then for $ 1\le j\le n$, we have,
$$
\varphi_{ij}=(\,(i-1)n+1+j-i,\ i\,n+1+j-i, \ldots
,(r-1)n+1+j-i,\ [rn+2-i,rn]\,)\,.
$$
% (In particular, we have,
%$$\varphi_{i1}=((i-1)n+2-i,i\,n+2-i, \cdots
%,(r-1)n+2-i,[rn+2-i,rn])$$
Here the succesive differences between the first $(r\sh+1\sh-i)$
entries is $n$, and the remaining entries are consecutive integers.
%last $(i\sh-1)$ entries are the interval of integers
%$[rn\sh+2\sh-i,\ rn]$.
%\section{A basis for $H^0(X,L_0)$:} In this section, we construct
%a natural basis for $H^0(X,L_0), X$ being any one of the $\tau$'s
%in the above chain.
%Towards proving the linear independence of admissible Pl\"ucker
%Coordinates on $X(w_s)$, we first prove some degree-$2$
%straightening relations among the admissible Pl\"ucker
%co-ordinates on $X(w_s)$.
\begin{lem}\label{dim}
\begin{enumerate}
\item $\# \{\, \text{\rm admissible $S$ on } \tau_1 \,\}=n^r$
\item For $1\le j\le n$,\ \
 $\# \{\, \text{\rm admissible $S$ on } \tau_j \,\}=n^{r-1}(n\sh+1\sh-j)$.
 \end{enumerate}
\end{lem}

\begin{proof} (1) Let $S:=(s_1, \cdots ,s_r)$
be a typical admissible tuple on $\tau_1=w_s$.
The expression for $w_s$ as a $r$-tuple
% (cf. Remark \ref{r-tuple})
together with the admissibility of $S$ implies
that there are $n$ ways of choosing $s_r$ from
$\{(r-1)n+1,(r-1)n+2,\cdots ,rn\}$; having chosen $s_r$, there are
$n$ ways of choosing $s_{r-1}$ from $\{s_{r}-n,\cdots ,s_{r}-1\}$,
and so on.
\\[.2em]
(2) The proof is similar: at the first step, we
have only $(n+1-j)$-choices for $s_r$ (since $s_r$ has to be $\ge
(r-1)n+j $, there are $(n+1-j)$ ways of choosing $s_r$ from
$\{(r-1)n+j,\cdots ,rn\}$). The rest is as before.
\end{proof}

\ni More generally, we have:

\begin{lem}
\begin{enumerate}
\item $\# \{\, \text{\rm admissible $S$ on } \varphi_{i1}\,\}=n^{r+1-i}$.
\item$\# \{\, \text{\rm admissible $S$ on } \varphi_{ij}\,\}=n^{r-i}(n\sh-j\sh+1)$.
\end{enumerate}
\end{lem}

\begin{proof}
(1) Since the last $(i-1)$ entries in $\varphi_{i1}$ are comprised
of $[rn+2-i,rn]$, we have that in any admissible $S$ on
$X(\varphi_{i1})$, the last $(i-1)$ entries are again comprised of
$[rn+2-i,rn]$. With $s_{r+1-i}$, we have $n$ choices (note that
$s_{r+1-i}$ may be chosen to be any entry from
\ni $[(r-1)n+2-i,rn+1-i]$). Rest of the proof is as in Lemma
\ref{dim}
\\[.2em]
(2) Again there is precisely one choice for $s_{r+2-i},\cdots
,s_{r}$. With $s_{r+1-i}$, we have $n+1-j$ choices (note that
$s_{r+1-i}$ may be chosen to be any entry from
$[(r-1)n+1+j-i,rn+1-i]$). The rest of the proof is as in Lemma
\ref{dim}
\end{proof}

\ni {\bf Pivots.} Recall that $1,\ n\sh+1,\ 2n\sh+1,\cdots
,(r\sh-2)n\sh+1,\ (r\sh-1)n\sh+1$, the entries of $w_s$,
are called the {\em pivots} of $w_s$.

For $1\le i\le r(n\sh-1)\sh+1$, let ${\mathcal{A}}_i$ denote
the set of all admissible $S:=(s_1, \cdots ,s_r)$ on $\tau_i$,
so that ${\mathcal{A}}_i\subset {\mathcal{A}}_{i-1}$.
We shall denote ${\mathcal{A}}_1$ also by just
${\mathcal{A}}$. If $S=(s_1, \cdots ,s_r) \in {\mathcal{A}}$
correponds to a Weyl group element $\tau $, then we shall denote
$s_k$ also by $\tau(k)$.
Let
$$
Z=\{S\in {\mathcal{A}}\mid s_r=(r\sh-1)n\sh+1 =\tau_1(r)\}$$
\begin{rem}\label{card} \begin{enumerate}
 \item \#$\,Z=n^{r-1}$

\item Let $S\in {\mathcal{A}}$.
 Then $p_{S}|_{X(\tau_{2})}$ is non-zero if and only
 if $S\not\in Z$. This is clear since:
 $$\tau_2=(2,\ n\sh+2,\ 2n\sh+2,\ldots ,(r\sh-2)n\sh+2,\ (r\sh-1)n\sh+2)\,,$$
  and hence  $p_{S}|_{X(\tau_{2})}$ is non-zero if and only
 if $s_r\ge (r-1)n+2$.  Note that $\tau_2$ is the smallest
 (under $\ge$) in ${\mathcal{A}}$ such that $s_r\ge (r\sh-1)n\sh+2$.

 \item ${\mathcal{A}}=Z\,\dot\cup\,{\mathcal{A}}_2$
 \end{enumerate}\end{rem}

\ni For $0\le j\le r-1$, set
 $$Z_j:=\left\{S=(s_1, \cdots ,s_r)\in Z\ \left|
 \begin{array}{c} s_i=(i\sh-1)n\sh+1=\tau_1(i)\\
 \text{for}\ j<i\le r\,,\\
 \text{and}\ s_j>(j\sh-1)n\sh+1
 \end{array}\right.\right\}$$
 i.e., for $S\in Z_j$,\ $j$ is the largest such that
 $s_j$ is not a pivot.   Note that:
 \begin{enumerate}
 \item $Z=\dot\cup_{0\le j\le
r-1}\,Z_j$.\\[-.9em]
 \item $Z_0=\{w_s\}$
% \ni 3. $Z_j$ is the set of $S\in{\mathcal{A}}$ with precisely
 %$n-j$ pivots.
 \end{enumerate}

 Let us arrange the elements of $Z$ lexicographically as
 $\{S_1<S_2<\cdots\}$, where
 $$
 S_1=(1,\, n\sh+1,\, 2n\sh+1,\ldots ,(r\sh-2)n\sh+1,\, (r\sh-1)n\sh+1)=\tau_1\,.
 $$
 $$
 S_2=(2,\, n\sh+1,\, 2n\sh+1,\ldots ,(r\sh-2)n\sh+1,\, (r\sh-1)n\sh+1)\,,
 $$
 etc,. The element  $S_2$ will play a crucial role in the discussion below.

\subsection{Degree two straightening relations}

In this subsection, we prove certain degree-two straightening relations
among admissible Pl\"ucker coordinates on $X(w_s)$ which will be
used for proving the linear independence of admissible Pl\"ucker
coordinates on $X(w_s)$ in the next subsection.

 We shall work with the open cell  $X^{\circ}(w_{s})$ in
 $X(w_{s})$ using the notation of \S\ref{ele},
 particulary the normalized $rn\times r$ matrix presentation $M$
 of a space $V\in X^\circ(w_s)$.
%Then as seen in \S \ref{ele}, $M=(a_{ij})$ may be thought of as
%(a lower triangular)
 % $rn\times r$ matrix, composed of $r$ blocks $A_1,\cdots ,A_r$ of $n\times r$
 % matrices such that the $1^{{\text{st}}}$ column of $A_1$ repeats as the $2^{{\text{nd}}}$
 % column of  $A_2$, repeats as the $3^{{\text{rd}}}$
 % column of  $A_3$, and so on.
 % \ni Similarly, the $1^{{\text{st}}}$ column of $A_2$ repeats as the $2^{{\text{nd}}}$
 % column of  $A_3$, repeats as the $3^{{\text{rd}}}$
 % column of  $A_4$, and so on.
 % \ni Further, we may suppose (by column
 % operations) that in the $((i-1)n+1)^{{\text{th}}}$ row (a pivotal
 % row), $a_{(i-1)n+1\,j}=0,\ j\not= i$, i.e.,
 %  $a_{(i-1)n+1\,i}$ is the only non-zero
 % entry.  With these conventions,
 We shall now compute
  $f_S:=p_{S}|_{X^{\circ}(w_{s}})$, for $S\in {\mathcal{A}}$ in terms
  of the entries in $M$. Further, if $S$ corresponds to a Weyl
  group element $\tau$, we shall denote $f_S$ also by $f_\tau$.
  %Let $x:=a_{11}$.

  \begin{lem}\label{special}
  We have:\ \
 (1)\ $f_{\tau_{1}}=1$ ;\quad %x^r$
 (2)\ $f_{\tau_{2}}=a_{21}^r$ ;\quad
 (3) $f_{S_{2}}=a_{21}$ %x^{r-1}$
  \end{lem}

  \begin{proof}
%The expressions for $f_{\tau_{1}},f_{\tau_{2}},f_{S_{2}}$ are
Clear by direct computations with the matrix $M$.
  \end{proof}

\begin{lem}\label{special2}
Let $S=(s_1, \cdots ,s_r)$ be in $Z_j$. Let $S'$ be obtained from
$S$ by replacing $s_{j+1}$ by $s_{j+1}+1$.
\begin{enumerate}
\item $S'$ is admissible; further, $S'$ is in $Z_{j+1}$ or
${\mathcal{A}}_2$ according as $j<r\sh-1$ or $j=r\sh-1$.
\item  $f_{S_{2}}f_{S}=f_{\tau_{1}}f_{S'}$
\end{enumerate}
\end{lem}
\begin{proof}
(1) The admissibility of $S'$ follows from the facts that $s_j$ is
not a pivot and $s_{j+1}$ is a pivot (since $S\in Z_j$); the
latter assertion in (1) is clear from this.
\\[.2em]
(2) We compare the evaluations of $f_S,f_{S'}$ on
the matrix $M$. We have $f_S(M)$ is the determinant of the matrix
$$\begin{pmatrix}M_{j}&0_{j,\, r-j}\\ *&D_{r-j,\,r-j}\end{pmatrix}$$
where $M_{j}$ is the $j\times j$ sub matrix of $M$ with row
indices $s_1,s_2,\cdots ,s_j$ and column indices $1,2,\cdots ,j$,
and $D_{r-j,\,r-j}$ is the diagonal $(r\sh-j)\times (r\sh-j)$ matrix
diag$\,(x,x,\cdots ,x)$. Hence
$$f_S(M)=x^{r-j}\Delta,\ \leqno{(*)}$$
where $\Delta:=\det M_{j}$.

Similarly, $f_{S'}(M)$ is the determinant of the matrix
$$\begin{pmatrix}M'_{j}&0_{j+1,\, r-j-1}\\
*&D_{r-j-1,\,r-j-1}\end{pmatrix}$$
where $M'_{j}$ is the $(j\sh+1)\times
(j\sh+1)$ submatrix of $M$ with row indices $s_1,s_2,\cdots
,s_j,s_{j+1}+1$ and column indices $1,2,\cdots ,j\sh+1$, and
$D_{r-j-1,\,r-j-1}$ is the diagonal $(r\sh-j\sh-1)\times (r\sh-j\sh-1)$ matrix
diag$\,(x,x,\cdots ,x)$. Hence
$$f_{S'}(M)=x^{r-j-1}\Delta',\ \leqno{(**)}$$
where $\Delta':=\det M'_{j}$.

But now the fact that $S,S'$ differ just in the $(j\sh+1)$-th entry implies that
the first $j$ columns of $M'_j$ are obtained by adding
$a_{nj+2,\,1}$, $a_{(n-1)j+2,\,1},\cdots $, $a_{n+2,\,1}$ respectively to
the $j$ columns of $M_j$
(in view of $\tau$-stability of $w_s$)
% (cf. Remark \ref{t-stability})),
and the last column of $M'$ consists
of zeroes except the last entry which is $a_{2,1}$.
 (Note that for
$1\le i\le r$, in the $i$-th column of $A$, the $((i\sh-1)n\sh+1)$-th
entry is $a_{11}=1$, and the $((i\sh-1)n\sh+2)$-th entry is $a_{21}$.)

We obtain $\Delta'=a_{21}\Delta$,
so that $(**)$ implies $f_{S'}(M)=x^{r-j-1}a_{21}\Delta$. Thus:
$$xf_{S'}(M)=x^{r-j}a_{21}\Delta=a_{21}f_S(M)\,,$$
which, with Lemma \ref{special}, implies:
$$f_{S_{2}}(M)f_{S}(M)=
a_{21}f_{S}(M)=f_{S'}(M)=f_{\tau_{1}}(M)f_{S'}(M)$$
From this (2) follows.
\end{proof}

More generally, we have similar quadratic relations among
$\{f_S,\,S\in {\mathcal{A}}_i\}$ for all $1\le i\le r(n\sh-1)$ as
given by the Lemma below. Let us first fix some notation. Fix a
$\varphi:=\tau_{l}$ for $1\le l\le r(n\sh-1)$. We can identify
$\varphi$ with some $\varphi_{ik}$ for $1\le i \le r$ ,
$1\le k\le n$ ; in fact, we may suppose that $k<n$, since
$\varphi_{in}=\varphi_{i+1,\,1}$ (since for $i=r$
we have $\varphi_{rn}=\mathrm{id}$).

Let us enumerate the elements of
${\mathcal{A}}_\varphi={\mathcal{A}}_l$ as $\{R_1,R_2,\ldots\}$
so that $\varphi=R_1<R_2<\cdots $ .
We have:
$$\varphi=\varphi_{ik}=
((i\sh-1)n\sh+1\sh+k\sh-i,\, in\sh+1\sh+k\sh-i, \ldots
,(r\sh-1)n\sh+1\sh+k\sh-i,\,[rn\sh+2\sh-i,rn] )\,.$$
Write $(i-1)n+1+k-i=p_{ik}n+q_{ik}$
where $1\le q_{ik}\le n$; for simplicity of notation, let us
denote $p_{ik},q_{ik}$ by just $p,q$ respectively. Then
$$\varphi=(pn+q, (p+1)n+q,\cdots,(r+p-i)n+q,[rn+2-i,rn])\,.$$

As with $\tau_1$, we shall do computations on the open
cell $X^{\circ}(\varphi)$ in
 $X(\varphi)$.
 The ${\mathbb{Z}}\times
 {\mathbb{Z}_+}$ matrix presentation for the elements of $X^{\circ}(\varphi)$
  has the following description:
 The ``relevant part" of the matrix presentation of a
 a generic point in $X^{\circ}(\varphi)$,\
 may be thought of as a lower triangular
  $(r-p)n\times r$ matrix $M:=(a_{gh})$
   composed of $r-p$ blocks $A_{p+1},\cdots ,A_{r}$
   of $n\times r$ matrices, where all the columns of $A_{p+1}$ except
   the first column consist of zeroes and the first column has the form
$$   \begin{pmatrix}
0\\
\vdots\\
a_{pn+q,\,1}\\
\vdots\\
a_{pn+n,\,1}
   \end{pmatrix}$$
Further, we have (in view of $\tau$-stability of $\varphi$)
that the $1^{{\text{st}}}$ column of $A_{p+1}$
repeats as the $2^{{\text{nd}}}$
  column of  $A_{p+2}$, repeats as the $3^{{\text{rd}}}$
  column of  $A_{p+3}$, and so on.
Similarly, the $1^{{\text{st}}}$ column of $A_{p+2}$ repeats as the $2^{{\text{nd}}}$
  column of  $A_{p+3}$, repeats as the $3^{{\text{rd}}}$
  column of  $A_{p+4}$, and so on.

Further, we define the pivotal rows to be those with indices
$\{mn+q\mid p\le m\le r\sh+p\sh-i\}$,
   together with all the entries of $[rn\sh+2\sh-i,rn]$.
In these rows, the matrix $M$ has a single entry 1 and all other
entries 0.
%Denoting $a_{pn+q,\,1}=x$,
That is:
$$\begin{gathered}
a_{mn+q,\,m+1-p}=1,\  \ p\le m\le r\sh+p\sh-i\\
  a_{rn+s-i,\,r+s-i}=1,\ \ 2\le s\le i
\end{gathered}$$
% Also since $x\not= 0$, we may suppose (by column
 % operations) that in a pivotal row, say,
%  $(mn+q)^{{\text{th}}}, a_{mn+q,m+1-p}$ is the only non-zero
 % entry. Thus
 $$a_{(mn+q)\,j}=0,\ \ j\not= m\sh+1\sh-p\,.$$
% (Corresponding to the pivotal row
%  $rn+s-i,\,a_{rn+s-i,\,r+s-i}$ is the only non-zero
%
  entry, $2\le s\le i$.)
In the discussion below, $f_S$ will denote $p_S\,|\,_{X^{\circ}(\varphi)}$.
Let $\varphi'$ denote $\varphi_{i\,k+1}$. Note that $X(\varphi')$ is a
divisor in $X(\varphi)$.

\begin{lem}\label{special3}
(1)\,$f_{\varphi}=1$\,;\ \ %x^{r+1-i}$
(2)\,$f_{\varphi'}=a_{pn+q+1,\,1}^{r+1-i}$\,;\ \
(3)\,$f_{R_{2}}=a_{pn+q+1,\,1}$\,.
%$f_{R_{2}}=x^{r-i}a_{pn+q+1,\,1}$
  \end{lem}

\ni The proof is similar to that of Lemma \ref{special}.

\begin{rem}\label{card2}
Let $$Z(\varphi)=\{S\in
{\mathcal{A}}_\varphi\mid s_{r+1-i}=\varphi(r\sh+1\sh-i)\}\,.$$
Then we have:
\\[.3em]
1. \#$\,Z(\varphi)=n^{r-i}$
 \\[.3em]
2. Let $S\in {\mathcal{A}}_\varphi$.
 Then $p_{S}|_{X(\varphi')}$ is non-zero if and only
 if $S\not\in Z(\varphi)$.
 This is clear because
 $\varphi'(r\sh+1\sh-i)=\varphi(r\sh+1\sh-i)+1$, and hence
 $p_{S}|_{X(\varphi')}$ is non-zero if and only
 if $s_{r+1-i}\ge \varphi(r\sh+1\sh-i)+1$.
 Indeed, $\varphi'$ is the smallest
 (under $\ge$) in ${\mathcal{A}}_\varphi$ such that
 $s_{r+1-i}\ge \varphi(r\sh+1\sh-i)\sh+1$.
 \\[.3em]
3. ${\mathcal{A}}_\varphi=Z(\varphi)\,\dot\cup\,{\mathcal{A}}_{\varphi'}$ .

 \end{rem}

\ni For $0\le j\le r\sh-i$, set:
 $$
 Z_j(\varphi):=\left\{S=(s_1, \cdots ,s_r)\in Z(\varphi)
 \ \left|\begin{array}{c}
  s_m=\varphi(m)\\ \text{for}\ j<m\le r\sh+1\sh-i,\\
  \text{and}\  s_j>\varphi(j)
  \end{array}\right.\right\}\,,
  $$
 i.e., for $S\in Z_j(\varphi)$,\ $j$ is the largest such that
 $s_j>\varphi(j)$. Note:
\\[.3em]
1. $Z(\varphi)=\dot\cup_{0\le j\le
r-i}\,Z_j(\varphi)$.
\\[.3em]
2. $Z_0=\{\varphi\}$

\begin{lem}\label{special4}
Let $S=(s_1, \cdots ,s_r)$ be in $Z_j(\varphi)$.
Let $S'$ be obtained
from $S$ by replacing $s_{j+1}$ by $s_{j+1}\sh+1$.
\begin{enumerate}
\item $S'$ is admissible; further, $S'$ is in $Z_{j+1}(\varphi)$ or
${\mathcal{A}}_{\varphi'}$ according as $j<r\sh-i$ or $j=r\sh-i$.
\item  $f_{R_{2}}f_{S}=f_{\varphi}f_{S'}$
\end{enumerate}
\end{lem}

\ni The proof is similar to that of Lemma \ref{special2}.

\subsection{Linear independence of admissible Pl\"ucker  coordinates}
\label{linind}

To show the independence of $\{p_S\mid S\in \mathcal{A}_\varphi \}$,
we first prove the independence of $\{p_S\mid S\in Z_j(\varphi) \}$
for $0\le j\le r\sh-i$.  We use induction on $\dim X(\varphi)$,
the starting point being $\varphi=\mathrm{id}$, for which the result is clear.

 Let then dim $X(\varphi)>0$.

 \ni First, let $j<r\sh-i$, and suppose
 $$
 \sum_{S\in Z_{j}(\varphi)}\ a_Sf_S=0,\ \ a_S\in  K\,. \leqno{(*)}
 $$
 Multiplying the relation by $f_{R_{2}}$,  by Lemma
 \ref{special4} it reduces to:
 $$
 f_{\varphi}\sum_{S'\in Z_{j+1}(\varphi)}\
 a_Sf_{S'}=0\,.
 $$
 $S'$ being as in Lemma
 \ref{special4}. Hence cancelling $f_\varphi$, it reduces to
  $$
  \sum_{S'\in Z_{j+1}(\varphi)}\  a_Sf_{S'}=0\,.
  $$
  Hence by decreasing induction on $j$, we obtain $a_S=0$, for all $S$.

  Let now $j=r-i$. Multiplying the relation ($*$)
  by $f_{R_{2}}$, it reduces to (in view of Lemma
 \ref{special4}):
  $$
 f_{\varphi}\sum_{S'\in{\mathcal{A}}_{\varphi'}} \  a_Sf_{S'}=0\,.
 $$
 Cancelling $f_\varphi$, and restricting to $X(\varphi')$,
 we obtain by induction on dim $X(\varphi)$ that
 $a_S=0$ for all $S$.

 Next, we prove the linear independence of $\{p_S\mid S\in Z(\varphi) \}$.
 Suppose:
 $$
 \sum_{S\in Z(\varphi)}\ a_Sf_S=0,\ \ a_S\in
 K\,. \leqno{(**)}$$
 As above, multiplying the relation (**) by
 $f_{R_{2}}$, cancelling $f_\varphi$ and restricting to $X(\varphi')$,
 we obtain $a_S=0$, for all $S\in Z_{r-i}(\varphi)$.
 Hence (**)  reduces to:
$$
\mathop{\sum_{S\in Z_{j}(\varphi)}}_{j<r-i}\
a_Sf_S=0,\ \ a_S\in K\,.\ \leqno{(***)}
$$
 Again multiplying the relation (***) by
 $f_{R_{2}}$, cancelling $f_\varphi$ and using the first step above,
 we obtain $a_S=0$, for all $S\in Z_{r-i-1}(\varphi)$.
 (Note that in the resulting relation, $\{a_S\mid S\in Z_{r-i-1}(\varphi)\}$
 occur as coefficients of the corresponding $S'\in Z_{r-i}(\varphi)$.)\
  Thus
 proceeding, we obtain $a_S=0$ for all $S\in
 Z(\varphi)$ appearing in (**).  Thus we obtain:

 \begin{prop}\label{linear}
$Z(\varphi)$ is linearly independent.
 \end{prop}

 We next prove the linear independence of
 $\{p_S\mid S\in{\mathcal{A}}_\varphi \}$.
 Since
 ${\mathcal{A}}_\varphi=Z(\varphi)\dot\cup\,{\mathcal{A}}_{\varphi'}$,
 we  may write a linear relation as:
 $$
 \sum_{S\in Z(\varphi)}\ a_Sf_S\,+
 \sum_{S\in {\mathcal{A}}_{\varphi'}}\ b_Sf_S=0\,.\leqno(\dagger)
 $$
 Restricting $(\dagger)$ to $X(\varphi')$, we first conclude that
 $b_S=0$, for all $S\in {\mathcal{A}}_{\varphi'}$.
 Then $(\dagger)$  reduces to:
 $$
 \sum_{S\in Z(\varphi)}\ a_Sf_S=0
 \,.$$
 In view of Proposition  \ref{linear}, it follows that $a_S=0$,
 for all $S\in Z(\varphi)$. Thus we obtain:

\begin{prop}\label{lin}
$\{p_S\mid S\in{\mathcal{A}}_{\tau_{\ell}} \}$ is linearly independent
for all $1\le \ell\le (r\sh-1)n\sh+1$.
 \end{prop}

 Proposition \ref{lin} together with the
generation by admissible $p_S$ (Proposition \ref{gen})
implies:

\begin{thm}\label{main1}
$\{p_S,S\in{\mathcal{A}}_{\tau_{\ell}}\}$ is a basis for
$H^0(X(\tau_\ell),L_0)$ for all $1\le \ell\le (r\sh-1)n\sh+1$.
\end{thm}

As an immediate Corollary, we obtain:
\begin{cor}\label{generation}
The straightening relations on $X(\tau_l)$ as given by Proposition
\ref{gen} give a set of generators for the the degree one part of
the ideal of $X(\tau_l)$ considered as a closed subvariety of
$Gr(sn,\,V_s)$.
\end{cor}

 As another consequence, we have:

 \begin{thm}\label{main2}
The shuffle relations  among the Pl\"ucker co-ordinates (of
Corollary \ref{shufflerels}) give a set of generators for the
degree one part of the ideal defining the affine Grassmannian
inside the infinite Grassmannian.
 \end{thm}

 \begin{proof}
As seen in \S\ref{ind} (Proposition \ref{pluckervanish},
Remark\ref{cano}), we have that the affine Grassmannian is the
inductive limit of the $X(w_s)$, and the result follows from
Theorem \ref{main1}, and Corollary \ref{generation}.
 \end{proof}

\ni{\textbf{Conjecture:}}\label{conj} The affine Grassmannian is
cut out inside the infinite Grassmannian by the (linear) shuffle
relations.

\section{Application to Nilpotent Orbit Closures}\label{sec_nilpotent}
\setcounter{subsubsection}{0}

\subsection{$\cp$-stable Affine Schubert Varieties}\label{stable}

Define

\noindent $W^P_{st}=\{w=(c_1,\ldots,c_n)\in W\mid
c_1+\cdots+c_n=0, c_1\leq\cdots\leq c_n\}$ $\subset W^P$.  One can
check that for $w\in W^P_{st}$, $X(w)$ is stable by left
translations by $\cp$ (and not just by $\cb$), and thus by
$\SL_n(K)\subset \cp$. Let $w\in W^P_{st}$. There exists a $s$
such that $X(w)\subset X(w_s)$ (in fact there are infinitely many
such $s$);  we consider $X(w)$ inside some fixed $X(w_s)$.

Define $\lam_w=(\lam_1,\ldots,\lam_n)\in \ZZ^n$ by
$\lam_i=-c_i+s$, $1\leq i\leq n$. Note that
$nc_n\geq\lam_1\geq\cdots\geq\lam_n\geq 0$ ($\lam_n\geq 0$ is
implied by $X(w)\subset X(w_s)$) and $\lam_1+\cdots+\lam_n=ns$.
Thus $\lam_w$ is a partition of $ns$ with at most $n$ non-zero
rows. Conversely, for any partition $\lam$ of $ns$ with at most
$n$ rows, we have that $\lam=\lam_w$ for a unique $w\in W^P_{st}$,
which we refer to as $w_{\lam}$ (we will often write $X(\lam)$ to
refer to $X(w_{\lam})$).  Thus there is a bijection between $\cp$
stable affine Schubert varieties in $X(w_s)$ and partitions of
$ns$ with at most $n$ rows, i.e., column lengths are $\le n$.

Let $\lam=(\lam_1,\ldots,\lam_m)$, $\mu=(\mu_1,\ldots,\mu_m)$ both
be partitions of $m\in\NN$. We say that $\lam\geq\mu$, if
$\sum\limits_{i=1}^{j}\lam_i\geq\sum\limits_{i=1}^{j}\mu_i$ for
$j=1,\ldots,m$. This order is referred to as the {\it dominance}
order. The following Lemma can be easily verified.

\begin{lem}\label{partitionorder}
Let $\lam$, $\mu$ be partitions of $sn$ with at most $n$ rows.
Then $w_{\lam}\geq w_{\mu}\iff \lam\geq\mu$ in the dominance
order.
\end{lem}

For $\mu$ a partition of $sn$ with at most $n$ rows, let $\lam$ be
the conjugate partition. Let $E=\CC^n$. The following result is
shown in \cite{Sh,mag2}.

\begin{thm}\label{shim}
$V_{d\Lambda_0}^{\star}(w_{\lambda})\cong
S_{\lambda_1^m}E\otimes\ldots\otimes S_{\lambda_s^m}E$ as $SL(E)$
modules (here, for a dominant weight $\nu=(\nu_1,\cdots,\nu_s),
L_\nu E$ denotes the associated Schur module).
\end{thm}

Applying Theorem \ref{th_irrrep_coordring}, we obtain the
following:

\begin{cor}\label{shimozono}
$K[X(\mu)]_d\cong L_{\lambda_1^m}E\otimes\ldots\otimes
L_{\lambda_s^m}E$ as $SL(E)$ modules.
\end{cor}

\subsection{Nilpotent orbit closures}

In the rest of this section, we discuss possible applications of
 our approach to nilpotent orbit closures in positive
 characteristics.

Let $\NNN$
 denote the set of all nilpotent matrices in $M_{n\times
n}(K)$; it is a closed affine subvariety of $M_{n\times n}(K)$.
The group $\GL_n(K)$ acts on $\NNN$ by conjugation. Each orbit
contains precisely one matrix in Jordan canonical form (up to
order of the Jordan blocks). Thus the orbits are indexed by
partitions $\mu$ of $n$, i.e. $\mu=(\mu_1,\ldots,\mu_n)$,
$n\geq\mu_1\geq\cdots\geq\mu_n\geq 0$, $\mu_1+\cdots+\mu_n=n$. The
orbit corresponding to the partition $\mu=(\mu_1,\ldots,\mu_n)$
will be denoted by $\NNN_{\mu}^0$, and its closure by
$\NNN_{\mu}$.
%Since $\NNN^0_{\lam}$ is irreducible, so is $\NNN_{\lambda}$.

Let $\lambda=(\lam_1,\ldots,\lam_n)$, the conjugate partition of
$\mu$. Then
$$
\NNN^0_{\mu}=\left\{N\in M_{n\times n}(K)\mid \hbox{
rank}(N^i)=n-\sum\limits_{j=1}^i\lam_j,i=1,\ldots,n\right\},
$$
\begin{equation*}\label{nilprankconds}
\NNN_{\mu}=\left\{N\in M_{n\times n}(K)\mid \hbox{ rank}(N^i)\leq
n-\sum\limits_{j=1}^i\lam_j,i=1,\ldots,n\right\}.
\end{equation*}
These two facts imply:

\begin{prop}\label{nilporbincl}
$\NNN_{\mu'}^0\subset \NNN_{\mu}\iff \mu'\leq\mu$ in the dominance
order.
\end{prop}

In particular, $\NNN_{(n)}$ contains all nilpotent orbits; thus it
equals $\NNN$.

Next we describe an isomorphism due to Lusztig between the nilpotent orbit closures
and open subsets of certain affine Schubert varieties.
Let $\kappa=w_s$ with $s=1$, so that $r=n$; thus, with the
convention in \S \ref{ele}, $\kappa$ corresponds to the $n$-tuple
$(1, n+1, 2n+1,\ldots, n(n-1)+1)$.
A generic point of $X(\kappa)$
has a $\ZZ\sh\times\ZZ_+$ matrix presentation with the relevant part
$M$ being a $n^2\times n$ matrix; further, $M$ consists of $n$
blocks of $n\times n$ matrices:
\begin{displaymath}
M=\left(\begin{array}{c} A_1\\ \vdots \\ A_{n-1}\\ A_n
\end{array}\right).
\end{displaymath}
Define $X'(\kappa)=\{M\in X(\kappa)\mid p_{L}(M)\neq 0\}$, where
$L$ is the $n$-tuple $(n(n-1)+1,\ldots,n^2)$ (note that
$p_{L}(M)=|A_n|$). Then $X'(\kappa)$ embeds into the affine
subspace of $M_{n^2\times n}(K)$ with lowest block equal to the
identity. For $\lam$ a partition of $n$, define
$X'(\lam)=X(\lam)\cap X'(\kappa)$ (here, $X(\lam)=X(w_{\lam})$).
%The varieties $X(\lam)$ all have the property that they are stable
%under left translation by $P_0$, and thus by $\GL_n(K)\subset
%P_0$. Explicitly, for $M\in X(\lam)$, $g\in \GL_n(K)$, g acts by
%left translation on each block of $M$.  This action restricts to
%an on $X'(\lam)$, on which the action is transitive.

Consider the morphism $\phi:\NNN\to M_{n^2\times n}$ given by
\begin{displaymath}
\phi(N)=\left(\begin{array}{c} N^{n-1}\\ \vdots \\ N \\ I
\end{array}\right).
\end{displaymath}
Lusztig has shown (cf. \cite{lu}, \cite{lu2}, \cite{mag})
%Then the image $\phi|_{\NNN_{\lam}}(\NNN_{\lam})$ is defined in
%$M_{n^2\times n}$ by the rank conditions in each block given by
%(\ref{nilprankconds}). One verifies that these are precisely the
%rank conditions defining $X'(\lam)$ as a subvariety of
%$M_{n^2\times n}(K)$.  Thus one has
%\begin{cor}
%$\phi:N_{\lam}\to X'(\lam)$ is an isomorphisms of affine
%varieties.
%\end{cor} One also verifies than $\phi$ is $\GL_n(K)$ equivariant.
%Let $n_{\lambda}$ be the matrix in lower Jordan canonical form
%with block sizes $\lam_1,\ldots,\lam_n$.
%\begin{lem}
%\begin{enumerate}
%\item $\phi(n_{\lambda})=e_{w_{\lam}}$, up to a permutation of the
%columns of $e_{w_{\lam}}$.
%\item $\phi$ is $\GL_n(K)$ equivariant.
%\item
%$\phi(\stab_{\GL_n(K)}(n_{\lam}))=\stab_{\GL_n(K)}(e_{w_{\lam}})$.
%\item $\overline{\GL_n(K)(n_{\lam})}=N_{\lam}$,
%$\overline{\GL_n(K)(e_{w_{\lam}})}= X'(\lam)$.
%\end{enumerate}
%\end{lem}

\begin{thm}\label{main3}
$\phi|_{\NNN_{\mu}}:\NNN_{\mu}\to X'(\mu)$ is an isomorphism of
affine varieties.
\end{thm}

For the rest of this section, we shall denote
$$E=K^n,\quad
Y=M_n(K)\,,\quad A= K[x_{i,j}]_{1\le i,j\le n}\,,$$ $K$ being our
algebraically closed field of arbitrary characteristic.  We
further denote $X:= (x_{i,j})$ the $n\times n$ generic matrix,
$J_\mu$ = the defining ideal of $\NNN_{\mu}$, and  $A_\mu
:=A/J_\mu$ the coordinate ring of $\NNN_{\mu}$.

We are interested in the properties of the coordinate ring
$A_\mu$ such as normality, rational singularities, equations, etc.
Since normality and rational singularities are known, we really
want to concentrate on the equations.

By the results of \cite{mehta} on the simultaneous Frobenius
splitting, the coordinate ring $A_\mu$ is normal and
Cohen-Macaulay (in fact the Springer type resolution of $\NNN_\mu$
is rational). Donkin (cf. \cite{donkin}) showed that the
coordinate ring of $\NNN$ has a good filtration i.e. it is
filtered by the Schur functors.

The defining equations are known (in a characteristic free way)
for the following two classes of partitions:
\\[.5em]
(a) $\mu = (r, 1^{n-r})$ is a hook. These are nilpotent matrices
of rank $\le r-1$. The defining ideal is generated by the $ r$
minors of a generic $n\times n$ matrix and the invariants (namely,
the coefficients of the characteristic polynomial). These
varieties are complete intersections in the corresponding
determinantal varieties.
\\[.5em]
(b) $\mu = (2^r ,1^{n-2r})$. These are the nilpotent matrices
with the square being zero, of rank $\le r$. In this case the
defining equations are the entries of the square of our matrix,
the invariants and the $ (r+1)$ minors. The minimal generators are
the entries of the square of the matrix and one irreducible
representation of highest weight $(1^{r+1} ,0^{n-2r-2},
(-1)^{r+1})$ given by cosets of $ (r+1)$ minors.
\\[-.5em]

Let $X$ be a generic $n\times n$ matrix. Consider the following
$n^2\times n$ matrix $M$:
%$$M= (I, X, X^2 ,\ldots ,X^{n-1}).$$
\begin{displaymath}
M:=M(X):=\left(\begin{array}{c}X^{n-1}\\ \vdots  \\X\\ I \\
\end{array}\right)
\end{displaymath}

 For a given $\mu$ we denote by ${\mathcal
F}_{m,\mu}$ the span in $A_\mu$ of the products of $\le m$ cosets
of maximal minors of the matrix $M$ in $A_\mu$. We have by
definition $${\mathcal F}_{m,\mu}\subset{\mathcal F}_{m+1,\mu}$$
We denote by $M(i_1 ,\ldots ,i_n)$ the maximal minor of $M$
corresponding to the rows $i_1 ,\ldots ,i_n$, for $1\le i_1
<\ldots <i_n\le n^2$.

Let us start by formulating a general conjecture regarding the
spaces ${\mathcal F}_{m,\mu}$. Let $E=K^n$.
\\[1em]
{\bf Conjecture:} Let $\lambda = (\lambda_1 ,\ldots ,\lambda_s
)$ be the partition conjugate to $\mu$. Then there is a
characteristic free isomorphism of $SL(E)$-modules
$$
{\mathcal F}_{m,\mu} =
L_{\lambda_1^m}E\otimes\ldots\otimes L_{\lambda_s^m}E
$$
where for a partition $\nu=(\nu_1 ,\ldots ,\nu_r ), L_\nu E$ denotes the
Weyl module with highest weight $(\nu_1 ,\ldots ,\nu_r )$.

\begin{rem}\label{char}
In characteristic $0$, the conjecture follows from Corollary
\ref{shimozono} (cf. \cite{mag2,Sh}).
\end{rem}

Our goal is to provide the explicit straightening of the products
of minors of $M$ to the ``standard products" of maximal minors.
Such straightening gives another presentation of the ring $A_\mu$
with the set of generators consisting of minors spanning
${\mathcal F}_{1,\mu}$ other than the minor $M(n(n-1)+1, n(n-1)+2,
\ldots ,n^2)$ which is equal to $1$. Notice that these generators
include the entries $x_{i,j}=\pm
M(n(n-2)+i,n(n-1)+1\ldots,n(n-1)+j-1,n(n-1)+j+1, \ldots ,n^2)$ for
$1\le i,j\le n$. In view of Theorem \ref{main3} we have:

\begin{prop} The straightening relations on $X(\mu)$ among the Pl\'cker
co-ordinates(with $M(n(n-1)+1, n(n-1)+2, \ldots ,n^2)$ specialized
to $1$) give a set of generators of the defining ideal $J_\mu$
(here, $X(\mu)$ is the affine Schubert variety associated to
$\mu$).
\end{prop}

 The trouble is that even if such a straightening is described, still
in order to get a set of generators of $J_\mu$ (in terms of the
matrix entries $x_{ij}$), we need to replace the generators of
${\mathcal F}_{1,\mu}$ other than the minors
$M(n(n-2)+i,n(n-1)+1\ldots,n(n-1)+j-1,n(n-1)+j+1, \ldots ,n^2)
(=x_{ij})$ for $1\le i,j\le n$, by suitable polynomial expressions
in the $x_{ij}$'s.

\subsection{Kostka-Foulkes polynomials}

We will describe the equations of the ideals $J_\mu$ over the field
of characteristic zero. This description was given in \cite{[W1]}.
Here we give a proof that brings out the role of the spaces
${\mathcal F}_{1,\mu }$.

Let us fix $n$ and $\mu$. Let us also fix a dominant integral
weight $\alpha = (\alpha_1 ,\ldots ,\alpha_n )$ for $GL(n)$,
namely, $\alpha_i\in \mathbb{Z}$ and $\alpha_1 \ge \alpha_2
\ge\ldots\ge\alpha_n$. We denote by $S_\alpha E$, the Schur module
associated to $\alpha$. The ring $A_\mu$ is a graded ring. Let us
denote by $m_{\alpha ,\mu ,i}$ the multiplicity of $S_\alpha E$ in
the $i$-th graded component of $A_\mu$. We define the series
$$
P_{\alpha ,\mu }(q):= \sum_{i\ge 0} m_{\alpha ,\mu
,i}\ q^i
$$

We call the series $P_{\alpha ,\mu} (q)$ the Poincare series
corresponding to the weight $\alpha$ and to the orbit
$\NNN_\mu^0$. Notice that the series $P_{\alpha ,\mu }(q)$ is non
zero only when $\alpha_1 +\ldots +\alpha_n =0$. We will assume
this throughout this section. All series $P_{\alpha ,\mu}(q)$ are
polynomials since by a formula of Kostant the total multiplicity
of $S_\alpha E$ in $A_{(n)}$ equals the multiplicity of the zero
weight in $S_\alpha E$.

The polynomials $P_{\alpha ,\mu }(q)$ have a very interesting
connection with the Kostka-Foulkes polynomials (see \cite{mac} for
the definition of Kostka-Foulkes polynomials). This connection is
related to Kraft's construction of intersecting with the diagonal
which we now describe.

Let $I\subset A (=K[x_{ij}])$ be the defining ideal of the set of
diagonal matrices, i.e. the ideal generated by the elements
$x_{i,j}$ with $i\ne j$. We define the algebras
$$B= A/I, B_{\mu'} := A/(I+J_\mu )$$
where $\mu^\prime$ is the partition conjugate to $\mu$; the reason
for using the conjugate partition $\mu^\prime$ to label the
algebra is that socle of $B_{\mu^\prime}$ is the Specht module
$\Sigma^{\mu^\prime}$).

  The algebras $B, B_{\mu'}$ are graded
algebras. There are  natural actions of the symmetric group $S_n$
on the graded algebras $B, B_{\mu'}$. Indeed, we can embed $S_n$
into ${GL(E)}$ by sending a permutation to the corresponding
permutation matrix. This defines an action of $S_n$ on $A$. The
ideals $I, I+J_\mu$ are $ S_n$-stable; hence we obtain natural
actions of $S_n$ on $B, B_{\mu'}$. These actions are compatible
with the natural surjections $B\rightarrow B_{\mu'}$,
$B_{\mu}\rightarrow B_{\nu}$ (where $\nu$ is less than $\mu$ in
the dominance order). We denote by $I_{\mu'}$ the kernel of the
surjection $B\rightarrow B_{\mu}$, i.e., $B_{\mu'}=B/I_{\mu'}$.

We have natural restriction maps $\psi :A\rightarrow B$, $\psi_\mu
:A_\mu\rightarrow B_{\mu'}$. We will denote the image of the
element $x_{i,i}$ in the algebras $B, B_{\mu'}$ by $y_i$. The
action of $S_n$ on $B, B_{\mu}$ is given by
$$\sigma (y_i )=y_{\sigma (i)}$$

For a partition $\lambda$ of $n$ we denote by $\Sigma^\lambda$ the
Specht module. Let $R$ be the representation ring of the symmetric
group $ S_n$. For a representation $V$ of $ S_n$ we denote by
$[V]$ its class in $R$. We also consider the polynomial ring
$R[q]$ over $R$ in a variable $q$. For a partition $\mu$ of $n$ we
define the element ${\tilde K}_\mu (q)\in R[q]$ (cf. \cite{mac})
as
$$ {\tilde K}_{\mu} (q)=\sum_{i\ge 0} [B_{\mu ,i}] q^i$$
where $B_{\mu ,i}$ denotes the $i$-th graded component of $B_\mu$.
We can write:
$$
 {\tilde K}_\mu (q)= \sum_{\lambda} {\tilde K}_{\lambda ,\mu} (q) [\Sigma^\lambda ]\,.
 $$

The polynomials ${\tilde K}_{\lambda ,\mu}(q)$ were extensively
studied in several contexts. In \cite{mac}, ch III. \S 6] they
were defined as the transition polynomials for the Hall-Littlewood
symmetric functions. In \cite{de-pro,[GP]} it was shown that the
definition above and the one in \cite{mac} agree. For two
partitions $\lambda,\mu$ denote $$a_{\mu ,\lambda ,i}:=\mathrm{\
the\ multiplicity\  of\ } \Sigma^\lambda\mathrm{\  in\  the\ }
i^{th}\mathrm{\  graded\  component\  of\  }B_{\mu}$$ Define
$$r_{\mu ,\lambda}(q)=\sum_{i\ge 0} a_{\mu ,\lambda ,i}q^i$$
 The following Proposition was conjectured in \cite{kra}
 and  proved in \cite{de-pro}.

\begin{prop}\label{symm}Let $\mu$ be a partition of $n$ and let $\mu^\prime
:=(\mu^\prime_1 ,\ldots ,\mu^\prime_t )$ be the conjugate
partition.
\begin{enumerate}\item The representation of $S_n$ on $B_\mu$ is
isomorphic to the induced representation from the trivial
representation of $S_{\mu^\prime_1}\times\ldots\times
S_{\mu^\prime_t}$ to $S_n$. \item{} The polynomial $r_{\mu
,\lambda}(q)$ is equal to the polynomial $\tilde{K}_{\mu
,\lambda}(q):= q^{-n(\lambda )}K_{\mu ,\lambda}(q^{-1})$ where
$K_{\mu ,\lambda}(q)$ is the Kostka-Foulkes polynomial(cf.
\cite{mac}).
\end{enumerate}
\end{prop}

The families of polynomials $P_{\alpha ,\mu}(q)$ and ${\tilde
K}_{\lambda ,\mu}(q)$ overlap in several important cases. The
first occurrence is when $\mu = (n)$. Let $\mu$ be a partition of
$n$ and let $\alpha = (\alpha_1 ,\ldots ,\alpha_n )$ be a dominant
weight for $GL(n)$ with $\alpha_1 +\ldots +\alpha_n =0$. Let
$\alpha_n =-k$. We can define a partition $\lambda :=\lambda
(\alpha )$ of $nk$ by setting $\lambda = (\alpha_1 +k,\ldots
,\alpha_n +k )$.

\begin{prop} (cf. \cite{[Gu],[He]}) We have
$P_{\alpha ,(n)}(q)={\tilde K}_{\lambda ,(k^n )}(q)$.
\end{prop}

The next occurrence deals with the following special kind of
dominant weights for $GL(E)$: let $\mathcal{W}_1^n$ denote the set
of dominant weights $\alpha = (\alpha_1 ,\ldots ,\alpha_n )$  for
$GL(E)$ with $\alpha_1 +\ldots +\alpha_n =0$, $\alpha_n =-1$. Such
weights are sometimes called {\it the weights of level one}. Let
$\mathcal{P}_n$ denote the partitions of $n$. We have a bijection
$$\diamond : \mathcal{W}_1^n\rightarrow \mathcal{P}_n,\ \diamond (\alpha ) = (\alpha_1
+1,\ldots ,\alpha_n +1 ) $$

\begin{thm}\label{level} Let $\mu$ be a partition of $n$ and let
$\alpha = (\alpha_1 ,\ldots ,\alpha_n )$ be in $\mathcal{W}_1^n$.
Let $\lambda =\diamond (\alpha )$, i.e. $\lambda = (\alpha_1
+1,\ldots ,\alpha_n +1 )$. Then we have
$$P_{\alpha ,\mu}(q)={\tilde K}_{\lambda ,\mu^\prime}(q)$$
\end{thm}

\begin{proof}  This follows from \cite{[W2]}, \S 6.
\end{proof}

Theorem \ref{level} gives another interpretation to the spaces
$\mathcal{F}_{m,\mu}$ related to the matrix $M$ defined at the
beginning of this section.

We now give a representation-theoretic interpretation for
${\mathcal F}_{1,\mu}$ in characteristic zero.

\begin{prop} Let $K$ be a field of characteristic zero. Let ${\mathcal
G}_{m,\mu}$ be the isotypic component of $A_\mu$ consisting of all
representations $S_{(\alpha_1 ,\ldots ,\alpha_n )}E$ with
$\alpha_n\ge -m$. Then ${\mathcal G}_{m,\mu}={\mathcal
F}_{m,\mu}$.
\end{prop}

\begin{proof} Let first $m=1$.
We know (by representation theory) that ${\mathcal
F}_{1,\mu}\subset {\mathcal G}_{1,\mu}$
%(note that ${\mathcalF}_{0,\mu}=K$).
Indeed, the maximal minors of the matrix $M$ span the subspace
whose $GL(n)$-representation type is a sub representation of
$\bigwedge^n E^*\otimes W$ where $W$ is a polynomial
representation of $GL(n)$ of degree $n$ (degree being with respect
to the matrix entries).

We will show that  for dimension reasons we in fact have equality.
For a weight $(\alpha_1 ,\ldots ,\alpha_n )$ with $\alpha_n =-1$,
$\alpha_1 +\ldots +\alpha_n =0$, consider the partition  $\nu
(\alpha ):= (\alpha_1 +1,\ldots ,\alpha_n +1) (=\diamond (\alpha
))$. Then by Theorem \ref{level},
%${\mathcal F}_{1,\mu}= {\mathcal G}_{1,\mu}$. Indeed,
the subspace ${\mathcal G}_{1,\mu}$ has the decomposition to
isotypic components given by Kostka-Foulkes polynomials:
$$
{\mathcal G}_{1,\mu} =
\mathop{\bigoplus_{\alpha\ :\ \alpha_n =-1}}_{ \alpha_1 +\cdots +\alpha_n
=0}\ \tilde{K}_{\mu' ,\nu (\alpha )}(1)\,S_{\alpha }E\,.
$$
Thus we see that as an $SL(E)$-module, we have
$$
{\mathcal G}_{1,\mu} = S_{\lambda_1}E\otimes\ldots\otimes S_{\lambda_s}E\,,
$$
$(\lambda_1,\cdots,\lambda_s)$ being the partition conjugate to
$\mu$; hence we obtain ${\mathcal G}_{1,\mu}={\mathcal
F}_{1,\mu}$. (Note that in view of Remark \ref{char} and Corollary
\ref{shimozono}, ${\mathcal
F}_{1,\mu}=L_{\lambda_1}E\otimes\ldots\otimes L_{\lambda_s}E$ .)\
\

For arbitrary $m$ it is clear that the subspace $\mathcal
{G}_{m,\mu}$ contains $\mathcal{ F}_{m,\mu}$. The equality is also
clear because the space $\mathcal{ F}_{1,\mu}$ generates $A_\mu$
as an algebra, and the equality is true for the case of the
nullcone $\mu = (n)$ (we have the dimension argument working for
arbitrary $m$ in that case).
\end{proof}

\begin{cor}\label{graded}
Let ${\mathcal F}_{{1,\mu},d} $ be the $d$-th graded component of
the space ${\mathcal F}_{1,\mu}$. we have
$$\sum_{d\ge 0} q^d \mathrm{char}({\mathcal F}_{{1,\mu},d} )=
\mathop{\bigoplus_{\alpha\ :\ \alpha_n =-1}}_{\alpha_1 +\cdots
+\alpha_n =0}\ \tilde{K}_{\mu' ,\nu (\alpha )}(q) S_{\alpha }E$$
\end{cor}

\subsection{Equations of nilpotent orbit closures}

Let us investigate the setting of Theorem \ref{level} more
closely. Notice that the Specht module $\Sigma^{\diamond (\alpha
)}$ can be defined as the set of vectors in $S_\alpha E$ of weight
zero. This means that the restriction map $\psi_\mu
:A_\mu\rightarrow B_{\mu'}$ takes the isotypic component of type
$\alpha$ in $A_\mu$ to the isotypic component of type $\lambda:=
\diamond (\alpha )$ in $B_{\mu'}$. Thus Theorem \ref{level} can be
strengthened as:

\begin{thm}\label{strong} Let $n, \lambda ,\alpha ,\mu $ be as
above. Let $A_\mu^\alpha$, $B_\mu^\lambda$ denote the isotypic
components of $A_\mu ,B_\mu$ respectively. Then the restriction
map induces  isomorphisms
$$\psi_\mu^\alpha :A_\mu^{\alpha}\rightarrow B_{\mu^\prime}^{\lambda}$$
\end{thm}

\begin{proof} First we prove the proposition for $\mu = (n)$. Recall that $T_1
,\ldots ,T_n$ (the coefficients of the characteristic polynomial
of the generic $n\times n$ matrix) are the basic invariants in
$A$. Denote $e_j :=e_j (y_1 ,\ldots ,y_n )$ the elementary
symmetric functions in $y_1 ,\ldots ,y_n$. Consider the Koszul
complexes $K(T_1 ,\ldots ,T_n ;A)$ and $K(e_1 ,\ldots ,e_n ; B)$.
The restriction map induces a map of complexes
$$\mathrm{res} : K(T_1 ,\ldots ,T_n ;A)^\alpha \rightarrow K(e_1 ,\ldots ,e_n ;B)^\lambda$$
which is an epimorphism on each term. Both complexes are acyclic
with zero homology groups being $A_{(n)}^\alpha$ and $B_{(1^n
)}^\lambda$ respectively. It follows by diagram chase that the
induced restriction map $\psi_{(n)}^\alpha :A_{(n)}^\alpha
\rightarrow B_{(1^n )}^\lambda$ is an epimorphism. Since the
dimensions of both spaces are the same by Theorem \ref{level}, the
map is an isomorphism as claimed.

To prove the general case we consider the commutative diagram
$$\begin{array}{clcr}
A_{(n )}^{\alpha }&
\stackrel{\psi^\alpha_{(n)}}{\rightarrow}&
B_{(1^n)}^{\lambda }\\
\downarrow&&\downarrow\\ A_{\mu}^{\alpha
}&\stackrel{\psi^\alpha_{\mu}}{\rightarrow} &B_{\mu^\prime
}^{\lambda } \end{array}$$ Since the vertical maps and the upper
horizontal map are epimorphisms it follows that
$\psi^\alpha_{\mu}$ is also an epimorphism. Again by Theorem
\ref{level} the vector spaces $A_{\mu}^{\alpha }$ and
$B_{\mu^\prime }^{\lambda }$ have the same dimension, so the map
$\psi^\alpha_{\mu}$ is an isomorphism and we are done.
\end{proof}

Let $A_\mu^d, B_{\mu^\prime}^d$ denote the elements of degree $d$
in $A_\mu, B_{\mu'}$ respectively. We will need the following
lemma regarding the multiplication by elements of degree one in
$A$ and $B$.

\begin{lem}\label{epi} Let $\mu, \alpha ,\lambda$, be as above.
Consider a representation $S_\alpha E$ contained in $A_\mu^d$ and
its restriction $\Sigma^\lambda$ contained in $B_{\mu^\prime}^d$.
Then the vector subspace $(\Sigma^\lambda ) B_{\mu^\prime}^1$ of
$B_{\mu^\prime}^{d+1}$ is the epimorphic image (by the restriction
maps in corresponding weights) of the vector subspace $(S_\alpha
E)A_\mu^1$ of $A_\mu^{d+1}$.
\end{lem}

\begin{proof} The vector space $B_{\mu^\prime}^1$ is
isomorphic to the Specht module $\Sigma^{(n-1,1)}$ (except for the
trivial case $\mu = (1^n )$). Let us investigate the tensor
product $\Sigma^\lambda \otimes_{K} \Sigma^{(n-1,1)}$ as a
$S_n$-module with the diagonal action of $S_n$. Let us denote by
$[\lambda ,1]$ the set of partitions which differ from $\lambda$
by exactly one box. Let $d(\lambda )$ denote the number of corner
boxes in $\lambda$ diminished by one. It is well known that the
above tensor product has the following decomposition.
$$\Sigma^\lambda \otimes_{ K} \Sigma^{(n-1,1)} =\oplus_{\nu\in [\lambda ,1]}\Sigma^\nu \oplus
(\Sigma^\lambda ) {}^{d(\lambda )} .$$

The vector space $A_\mu^1$ is isomorphic to $S_{(1,0^{n-2},-1)}E$
(except for the trivial case $\mu = (1^n )$).
Littlewood-Richardson Rule implies the tensor product
decomposition
$$
S_\alpha E\otimes S_{(1,0^{n-2},-1)}E
 = \oplus_{\nu\in [\lambda ,1]} S_{\diamond^{-1} (\nu )} E\oplus
(S_\alpha E)^{d(\lambda )}\oplus \oplus_{\beta\notin
\mathcal{P}(1)} S_\beta E^{\oplus m(\alpha ,\beta )} \,.
$$
Here $m(\alpha ,\beta )$ denotes the multiplicity of $S_\beta E$.

Moreover, since the Specht module $\Sigma^{\diamond (\alpha )}$ is
obtained from $S_\alpha E$ by restricting to weight $0$, one
concludes that the restriction to weight zero induces the
epimorphism from the first two summands in the second
decomposition to the tensor product $\Sigma^\lambda \otimes_{W}
\Sigma^{(n-1,1)}$. Let us call the sum of the first two summands
in the second decomposition by $(S_{\alpha} E\otimes
S_{(1,0^{n-2},-1)}E)_{\mathcal{P}(1)}$. The lemma follows by
a chase of the diagram:
%\[
%\begin{smallmatrix}\hspace{-1em}
%(S_\alpha E\otimes S_{(1,0^{n-2},-1)}E)_{{\mathcal
%P}(1)}&\hookrightarrow& S_\alpha E\otimes
%S_{(1,0^{n-2},-1)}E
%&\rightarrow&
%A_\mu^1\otimes A_\mu^d&
%\rightarrow&
%A_\mu^{d+1}\\[0em]
%&h\searrow&\downarrow&&\downarrow&&\downarrow\\[0em]
%&&
%{\Sigma}^{(n-1,1)}\otimes {\Sigma}^\lambda &\rightarrow&
%B_{\mu'}^1\otimes B_{\mu'}^d&\rightarrow&B_{\mu'}^{d+1}
%\end{smallmatrix}
%\]
%
$$\hspace{-1em}
\begin{array}{c@{\,}c@{\,}c@{\,}c@{\,}c@{\,}c@{\,}c}
(S_\alpha E\otimes S_{(1,0^{n-2},-1)}E)_{{\mathcal P}(1)}
&\hookrightarrow&
S_\alpha E\otimes S_{(1,0^{n-2},-1)}E&
\rightarrow&
 A_\mu^1\otimes A_\mu^d&
\rightarrow& A_\mu^{d+1}\\[.5em]
&{}_h\!\!\!\!\searrow&\downarrow&&\downarrow&&\downarrow\\[.5em]
&& {\Sigma}^{(n-1,1)}\otimes {\Sigma}^\lambda &\rightarrow&
B_{\mu'}^1\otimes B_{\mu'}^d&\rightarrow&B_{\mu'}^{d+1}
\end{array}$$
taking into account that the map $h$ and all vertical maps are
epimorphisms.
\end{proof}

Theorem \ref{strong} and Lemma \ref{epi} allow us to describe
explicitly the generators of the ideals $J_\mu$,
 thus giving a new proof of the result of the last author (cf. \cite{[W1]}).
 The main idea is that we know for other reasons that the equations
of $J_\mu$ have to belong to the subspace $\mathcal{F}_{1,\mu}$
and we know the generators of the ideals $I_{\mu'}$.

Let us define the spaces $V_{i,p}$. They are defined for $0\le
i\le \min (p,n\sh-p)$. The space $V_{i,p}$ is  contained in the
span of $p\times p$ minors of $X$ (which we identify with
$\wedge^p E\otimes\wedge^p E^*$). The subspace $V_{i,p}$ is
spanned by the image of the map
$$\wedge^i E\otimes\wedge^i E^*
\stackrel{1\otimes tr^{(p-i)}}{\longrightarrow} \wedge^i
E\otimes\wedge^i E^*\otimes\wedge^{p-i} E\otimes\wedge^{p-i}E^*
\longrightarrow\wedge^p E\otimes\wedge^p E^*$$ where the last map
is a tensor product of exterior multiplication. The space
$V_{i,p}$ can be alternatively defined as a span of linear
combinations
$$\sum_{|J|=p-i}  X(P, J|Q, J)$$
for all subsets $P, Q$ of $[1,n]$ of cardinality $i$. Here $X(P|
Q)$ denotes the minor of $X$ corresponding to rows indexed by $P$
and columns indexed by $Q$.

The main result of \cite{[W1]} is:

\begin{thm}\label{old} The ideal $J_\mu$ is generated by the spaces
$V_{0,p}$ ($1\le p\le n$) and by the spaces $V_{i,\mu (i)}$ ($1\le
i\le n$) where $\mu (i):=\mu_1 +\cdots +\mu_i -i+1$.
\end{thm}

Let us look at the equations defining the algebra
$B_{\mu^\prime}$. The restrictions of the representations
$V_{i,p}$ that vanish on $Y_\mu$ certainly map to zero in
$B_{\mu^\prime}$. The main point of our proof of Theorem \ref{old}
is that one can prove using
 only combinatorial means to show
 that the restrictions of the representations $V_{i,p}$
described in Theorem \ref{old} generate the ideal
$I_{\mu^\prime}$. More precisely, we will use the generators of
$I_{\mu^\prime}$ of (cf. \cite{[GP]}) described below in
Proposition \ref{ideal} (note that our $\mu$ is $\mu^\prime$ in
\cite{[GP]}).

Let $S=\lbrace j_1 ,\ldots ,j_k\rbrace$ be a subset of $[1,n]$.
For every $r$, $1\le r\le \# S$, we define $e_r (S)$ to be the
$r$-th elementary symmetric function in the variables $y_{j_{1}}
,\ldots ,y_{j_{k}}$.

\begin{prop}\label{ideal} The ideal $I_{\mu^\prime}$ is generated by the following set
$\mathcal{C}_\mu$ of symmetric functions on subsets of variables
$y_1 ,\ldots ,y_n$.
$$\mathcal{C}_\mu =\lbrace e_r (S)\mid k\ge r>k\sh-d_k (\mu ),\   \# S=k ,\  S\subset [1,n]\rbrace$$
where $d_k (\mu )=\mu_{n-k+1} +\cdots +\mu_n $.
\end{prop}

\begin{prop}\label{alg}The algebra $B_{\mu^\prime}$ is the factor of $B$ by the ideal
$I_{\mu^\prime}$ generated by the restrictions of the
representations $V_{i,\mu (i)}$ where $\mu (i)=\mu_1 +\ldots
+\mu_i -i+1$ (which are zero if $i>min (\mu (i), n-\mu (i) )$) and
by the elementary symmetric functions $e_j (y_1 ,\ldots ,y_n )$.
\end{prop}

\begin{proof} We compare the restrictions of
representations $V_{i,p}$ to the generators exhibited in
Proposition \ref{ideal}. The restriction of $V_{i,p}$ is the span
of elements of the form $y_{j_{1}}\ldots y_{j_{i}} e_{p-i}(
y_{j_{i+1}} ,\ldots ,y_{j_{p}} )$. where $j_1 ,\ldots ,j_p$ are
distinct elements of $[1,n]$.

Let us perform induction on $i$.
For $i=0$, we get all elementary symmetric functions in $y_1
,\ldots ,y_n$ which can generate all elements of $\mathcal{C}_\mu$
with $\# S=n$. For $i=1$, we get all the symmetrizations of the
functions $y_1 e_{p-1} (y_2 ,\ldots ,y_n )$. Modulo elementary
symmetric functions we can, however, write
$$y_1 e_{p-1} (y_2 ,\ldots ,y_n )\equiv -e_p (y_2 ,\ldots ,y_n )$$
because
$$e_p (y_1 ,\ldots ,y_n )=y_1 e_{p-1} (y_2 ,\ldots ,y_n )\equiv e_p (y_2 ,\ldots ,y_n ).$$
This gives us the functions in $\mathcal{C}_\mu$ corresponding to
$\# S=n-1$, since the condition $p>n-1-d_{n-1} (\mu )$ means
exactly $p>\mu_1 -1$.

\noindent For $i=2$ we get the symmetrizations of the functions
\noindent $y_1 y_2 e_{p-2}(y_3 ,\ldots ,y_n )$. But modulo the
elements which are restrictions of those from $V_{1,p}$ and
$V_{0,p}$ we have
$$y_1 y_2 e_{p-2}(y_3 ,\ldots ,y_n )\equiv -e_p (y_3 ,\ldots ,y_n )$$
because $e_p (y_1 ,\ldots ,y_n )$ can be written as
$$e_p (y_3 ,\ldots ,y_n )+y_1 e_{p-1}(y_3 ,\ldots ,y_n )+y_2 e_{p-1}(y_3 ,\ldots ,y_n )
+y_1 y_2 e_{p-2}(y_3,\ldots ,y_n )$$
 The condition
$r>n-2-d_{n-2}(\mu )$ means exactly $p>\mu_1 +\mu_2 -2$ so we
generate in this way all elements in $\mathcal{C}_\mu$ with $\#
S=n-2$.

Continuing in this way we prove that $\mathcal{C}_\mu$ is
contained in the ideal generated by the restrictions of the
representations $V_{i,p}$ listed in (8.2.5) of \cite{[W2]} which
proves the proposition.
\end{proof}

\noindent{\sl  Proof of Theorem \ref{old}}. First of all we need
to know that the geometric method of calculating syzygies applied
to a Springer type desingularization of $\NNN_\mu$ allows to show
that the ideal $J_\mu$ is generated by the representations of type
$S_{(1^i ,0^{n-2i}, (-1)^i )}E$ for $0\le i\le \frac{n}{2}$. This
is Corollary (8.1.7) in \cite{[W2]}.

We investigate the restriction map $\psi_\mu$ on the isotypic
spaces corresponding to weights $\alpha (i)=(1^i ,0^{n-2i}, (-1)^i
)$, $\lambda (i)= (2^i , 1^{n-2i})$. Consider the restriction map
$$\psi_\mu^{\alpha (i)}
:A_\mu^{\alpha (i)}\rightarrow B_{\mu^\prime}^{\lambda (i)}$$

This map is an isomorphism. Consider the diagram
$$\begin{array}{clcr}
A_{(n )}^{\alpha (i)}&
\stackrel{\psi_{(n)}^{\alpha (i)}}{\longrightarrow}
 & B_{(1^n )}^{\lambda (i)}\\
\downarrow&&\downarrow\\ A_{\mu}^{\alpha (i)}&
\stackrel{\psi_\mu^{\alpha (i)}}{\longrightarrow}&
B_{\mu^\prime}^{\lambda (i)}
\end{array}$$
where the horizontal maps are the
restriction isomorphisms and the vertical maps are natural
surjections.

Denote by $J_\mu^\prime$ the ideal generated by the
representations defined in Theorem \ref{old}.  By \cite{[W2]}, Lemma 8.2.1,
 we have that $J_\mu^\prime\subset J_\mu$.

Assume that there is a representation $S_{\alpha (i)} E$ in $A_{(n
)}^{\alpha (i)}$ which is not in the ideal $J^\prime_\mu$. Take
such a representation of minimal possible degree $d$. It is
therefore among the generators of $J_\mu$. It maps down to zero in
$A_{\mu}^{\alpha (i)}$ and therefore its image $\Sigma^{\lambda
(i)}$ in $B_{(n)}^{\lambda (i)}$ is generated by elements of lower
degree in $I_{\mu^\prime}$. Therefore there is a representation
$\Sigma^\nu$ in degree $d-1$ in $I_{\mu^\prime}$ for which
$\Sigma^{\lambda (i)}$ is contained in the product $(\Sigma^\nu
)S_{\mu^\prime}^1$.

This however means by Lemma \ref{epi} that $S_{\alpha (i)} E$ is
in the product $(S_\beta E) A_\mu^1$ where $S_\beta E$ is the
representation restricting to $\Sigma^\nu$. Since $S_\beta E$ is
in $J_\mu$, we get a representation in lower degree in $J_\mu$ but
not $J_\mu^\prime$. Since $J_\mu$ is generated by representations
of weights $\lambda (i)$, we can get in lower degree a
representation of a weight $\lambda (j)$ which is in $J_\mu$ but
not $J_\mu^\prime$. This gives a contradiction to the minimality
of $d$.

\begin{rem}
There are fascinating combinatorial expressions  for the
polynomials ${\tilde K}_{\lambda,\mu}(q)$ due to Lascoux and
Sch\"utzenberger (cf. \cite{[LS]}). They define a statistic, the
{\it charge}, on the
set $\mathrm{ST}(\lambda )_{\mu}$ of standard tableaux of shape $\lambda$
and of content $\mu$ (i.e. containing $\mu_1$ 1's, $\mu_2$ 2's,
and so on) such that:
$$
{\tilde K}_{\lambda ,\mu}(q)
=\sum_{T\in \mathrm{ST}(\lambda )_{\mu}} q^{\mathrm{charge}(T)}
$$
  A  general conjecture, giving the combinatorial description of
  $P_{\alpha ,\mu}(q)$ for arbitrary $\alpha$, is
described in \cite{[SW]} and proved for special partitions $\mu$.
\end{rem}

\end{document}